\documentclass[12pt]{amsart}   
\usepackage{amssymb,amsmath,amscd,verbatim}  
\usepackage[all]{xy}  
\font\tenrsf=rsfs10 at 11pt  
\font\sevenrsf=rsfs7 at 8pt  
\font\fiversf=rsfs5 at 6pt  
\newfam\rsffam  
\textfont\rsffam=\tenrsf  
\scriptfont\rsffam=\sevenrsf  
\scriptscriptfont\rsffam=\fiversf  
\def\rond#1{{\tenrsf\fam\rsffam#1}}  
\def\thesection{\arabic{section}} 
\renewcommand{\theequation}{\arabic{section}.\arabic{equation}} 
\newtheorem{theorem}{Theorem}[section] 
\newtheorem{lemma}[theorem]{Lemma} 
\newtheorem{proposition}[theorem]{Proposition} 
\newtheorem{corollary}[theorem]{Corollary} 
 
\theoremstyle{definition}  
\newtheorem{definition}[theorem]{Definition} 
\newtheorem{example}[theorem]{Example} 
\theoremstyle{remark}  
\newtheorem{rem}[theorem]{Remark}  

\newcommand{\Hr}{\rond{H}}

\newcommand{\Cc}{\mathcal{C}}  
\newcommand{\Dc}{\mathcal{D}}  
  
\newcommand{\Ic}{\mathcal{I}}  
\newcommand{\Jc}{\mathcal{J}}

\newcommand{\Sc}{\mathcal{S}} 
 
\newcommand{\Vc}{\mathcal{V}}  
\newcommand{\ad}{\mathrm{ad}}

\newcommand{\re}{\mathrm{Re}}   
\newcommand{\im}{\mathrm{Im}}   
\newcommand{\id}{\mathrm{Id}}

\newcommand{\supp}{\mathrm{supp}}  
\newcommand{\C}{\mathbb{C}}  
\newcommand{\R}{\mathbb{R}}  
\newcommand{\N}{\mathbb{N}}  
\newcommand{\Z}{\mathbb{Z}}

\def\build#1_#2^#3{\mathrel{\mathop{\kern 0pt#1}\limits_{#2}^{#3}}} 
\def\cchi{\raisebox{.45 ex}{$\chi$}}
\def\l{<}   
\def\g{>}  
\newcommand{\ind}{1} 

\begin{document}  
\title[A new look at Mourre's commutator theory]  
{A new look at Mourre's commutator theory}   
\author{Sylvain Gol\'enia}   
\address{Departement of Mathematical Methods in Physics\\ 
Warsaw University\\ Ho\.za 74\\ 
00-682 Warszawa, Poland} 
\email{Sylvain.Golenia@fuw.edu.pl}   
 \author{Thierry Jecko}   
\address{IRMAR, UMR 6625 du CNRS, Universit\'e de Rennes 1, Campus de Beaulieu,  
35042, Rennes Cedex France}   
\email{Thierry.Jecko@univ-rennes1.fr}   
\date{\today}   
\begin{abstract}   
Mourre's commutator theory is a powerful tool to study the continuous 
spectrum of self-adjoint operators and to develop scattering theory. We propose a 
new approach of its main result, namely the derivation of the limiting absorption principle 
(LAP) from a so called Mourre estimate. We provide a new interpretation 
of this result. 
\end{abstract}   
\maketitle   
\tableofcontents 
\section{Introduction} \label{intro}

In the beginning of the eigthies, Mourre's commutator theory was developed in \cite{m} 
to show absolute continuity of the continuous spectrum of $N$-body Schr\"odinger operators
and to study their scattering theory (cf. \cite{abg,hus}). In particular, one wanted to show 
their asymptotic completeness and the Mourre estimate (cf. \eqref{eq:mourre2}) played a crucial 
role in the proof (cf. \cite{dg,hus}). Now, Mourre's commutator theory is fundamental tool to develop the 
stationary scattering theory of general self-adjoint operators. We refer to \cite{abg,dg} for 
details. We point out that the theory is still used (see \cite{bchm, cgh, dj, ggo}, for instance) and 
that there were new developements to apply it to quantum field theory (cf. \cite{ggm1,ggm2}). 
The theory uses a so called {\em differential inequality technics}, that is quite magic and 
mysterious (to us at least). In this paper, we propose a new approach and interpretation of the 
theory. Since the original method has been developed to a rather sophisticated level (cf. 
\cite{abg,ggm1,s}), we did not try to optimize our approach and to give new results, but to 
focus on an intermediate, interesting situation. However, Theorem~\ref{th:mourre2} gives an 
extension of results in \cite{c,cgh}. We point out that our new approach of Mourre's 
commutator theory is an adaptation of a strategy to get semiclassical resolvent estimates for 
Schr\"odinger operators. This strategy was introduced by the second author in \cite{jec} and 
further used in \cite{cj,jec2}. 

To enter into the details of our approach, we need some notation and basic notions (see 
Subsection~\ref{base} for details). We consider two self-adjoint 
(unbounded) operators $H$ and $A$ acting in some complex Hilbert space $\Hr$. Let $\|\cdot\|$ denote 
the norm of bounded operators on $\Hr$. We shall study spectral properties of $H$ with the help 
of $A$. Since the commutator $[H, iA]$ is going to play a central role in the theory, 
we need some regularity of $H$ with respect to $A$ to give an appropriate sense to this commutator. 
Since $H$ is self-adjoint, its spectrum is included in $\R$. Given $k\in \N$, we say that $H\in 
\Cc^k(A)$ if for some (and thus for all) $z\not\in \R$, for all $f\in \Hr$, 
the map $\R\ni t\mapsto e^{itA}(H-z)^{-1}e^{-itA}f\in \Hr$ has the usual $\Cc^k$ regularity.
Let $H\in \Cc^1(A)$ and $\Ic$ be a bounded interval of $\R$. 
We say that the \emph{Mourre estimate} holds true for 
$H$ on $\Ic$ if there exist $c\g 0$ and a compact operator $K$ such that   
\begin{eqnarray}\label{eq:mourre2} 
E_\Ic (H)[H, iA]  E_\Ic(H)\geq c E_\Ic(H) + K, 
\end{eqnarray}  
in the form sense on $\Hr\times \Hr$. Here $E_\Ic(H)$ denotes the spectral measure of 
$H$ above $\Ic$.  
\begin{rem}\label{r:viriel+mourre}
Let $f\in \Hr$ and $\lambda\in \Ic$ with $Hf=\lambda f$. Then $E_\Ic(H)f=f$. Assume that $H\in\Cc^1(A)$. The Virial theorem (cf. \cite[Proposition 7.2.10]{abg}) implies that 
$\langle f,[H, iA]f\rangle =0$. If \eqref{eq:mourre2} holds true with $K=0$ then $f$ must be 
zero and there is no eigenvalue in $\Ic$. If \eqref{eq:mourre2} holds true then the total 
multiplicity of the eigenvalues in $\Ic$ is finite (cf. \cite[Corollary 7.2.11]{abg}). A weaker 
version of this result is due to Mourre in \cite{m}. For a general discussion on the Virial theorem 
see \cite{gg}. 
\end{rem}
The main aim of Mourre's commutator theory is to show the \emph{limiting absorption principle} 
(LAP) on some bounded interval $\Ic$ in $\R$. Given such a $\Ic$ and $s\geq 0$, we say that the   
LAP, respectively to the triplet $(\Ic,s,A)$, holds true for $H$ if  
\begin{equation}\label{eq:lap}  
\sup_{{\rm Re}z\in \Ic , {\rm Im}z\neq 0}\|\langle A \rangle^{-s}(H-z)^{-1} \langle A\rangle^{-s}\|<\infty.    
\end{equation}  

\begin{theorem}\label{th:mourre} 
Let $H\in \Cc^2(A)$, $\Ic$ be a bounded, open interval, and $s\g 1/2$. 
Assume the strict Mourre estimate, i. e. \eqref{eq:mourre2} with $K=0$, holds true. 
Then, for any closed subinterval $\Ic'$ of $\Ic$, the LAP for $H$ respectively to $(\Ic', s, A)$ holds true. 
\end{theorem}

\begin{rem}\label{r:K=0}
Assume the Mourre estimate \eqref{eq:mourre2} holds true on $\Ic$ with $K\neq 0$. Then, 
on small enough intervals outside the point spectrum $\sigma _{pp}(H)$ of $H$, which is finite 
by Remark~\ref{r:viriel+mourre}, the strict Mourre estimate \eqref{eq:mourre2} with $K=0$ 
holds true (cf. \cite{abg}) and Theorem~\ref{th:mourre} applies there. Putting all together, this 
yields the LAP on any compact subset of $\Ic\setminus\sigma _{pp}(H)$. 
\end{rem}

Compared with previous results, we do not need that the domain $\Dc(H)$ of $H$ is invariant under 
the $C_0$-group generated by $A$ (i.e. the propagator of $A$) or that $H$ has a spectral gap 
(cf. \cite{abg}). The main reason for this comes from the fact that we do not work with $H$ itself 
but with a local version of $H$, which is a bounded operator. This explains also why we can replace 
the global regularity assumption $H\in \Cc^2(A)$ by a local one and get a stronger result, 
namely Theorem~\ref{th:main}. The latter is covered by Sahbani's result in \cite{s} 
(cf. Remark~\ref{r:s}). Motivations for Theorem~\ref{th:main} are given in 
Remarks~\ref{r:s} and~\ref{r:exemple}. In Subsection~\ref{s:sketch}, we give a sketch of the 
proof of Theorem~\ref{th:main} and present our interpretation of Mourre's commutator theory, 
which is close to the interpretation of Remark~\ref{r:viriel+mourre}. We do not use the usual differential
inequality technics. 

In some sense, Theorem~\ref{th:mourre} (and also Theorem~\ref{th:main}) is not 
satisfactory (cf. Subsection~\ref{s:motiv}) and one wishes to replace the 
resolvent $(H-z)^{-1}$ in \eqref{eq:lap} by the reduced resolvent, namely 
$(H-z)^{-1}P^\perp$, where $P^\perp =1-P$, and $P$ is the orthogonal projection onto the 
pure point spectral subspace of $H$. 
For $s\geq 0$, we say that the {\em reduced} LAP, respectively to the triplet $(\Ic,s,A)$, holds true 
for $H$ if  
\begin{equation}\label{eq:tal-reduit}  
 \sup_{{\rm Re}z\in \Ic , {\rm Im}z\neq 0}\|\langle A \rangle^{-s}(H-z)^{-1} P^\perp
\langle A\rangle^{-s}\|<\infty .
\end{equation} 
\begin{theorem}\label{th:mourre2} 
Let $H\in \Cc^2(A)$, $\Ic$ be a bounded, open interval and let $s\g 1/2$. Assume the Mourre 
estimate \eqref{eq:mourre2} holds true on $\Ic$. Assume also that the range ${\rm Ran}PE_\Ic$ 
of $PE_\Ic$ is included in the domain $\Dc (A^2)$ of $A^2$. Then, for all closed interval $\Ic'$ 
included in the interior of $\Ic$, the reduced LAP \eqref{eq:tal-reduit}, respectively to $(\Ic ',s,A)$, 
holds true for $H$. 
\end{theorem} 

A similar result appears in \cite{cgh}. The authors assume a stronger regularity 
(essentially like $H\in \Cc^4(A)$) that implies ${\rm Ran}PE_\Ic\subset \Dc (A^2)$, by \cite{c}, 
and then show \eqref{eq:tal-reduit}.
The latter result and Theorem~\ref{th:mourre2} actually work with weaker, ``local" 
assumptions as shown in Proposition~\ref{p:mourre-regu}  and Theorem~\ref{th:mourre2'}. 
As mentioned before, this local version of the result might be important (cf. Remarks~\ref{r:s} and~\ref{r:exemple}). In \cite{cgh}, some H\"older continuity of the boundary values of the 
reduced resolvent $\lim _{\epsilon\to 0^+}\langle A \rangle^{-s}(H-\lambda -i\epsilon)^{-1} P^\perp\langle A\rangle^{-s}$ is obtained. In Remark~\ref{r:s-reduit} we explain how to get 
this under ``local" assumptions, using \cite{s}. \\
We point out that our proofs of Theorems~\ref{th:mourre2} and~\ref{th:mourre2'} is a quite 
immediate generalization of our proofs of Theorems~\ref{th:mourre} and~\ref{th:main}. We 
also give an alternative proof of Theorem~\ref{th:mourre2'} which is close to the corresponding 
proof in \cite{cgh}. Notice further that, Theorems~\ref{th:mourre2'} works under a {\em projected 
Mourre estimate} \eqref{eq:mourremain-proj}, that is weaker than \eqref{eq:mourre2}. In 
Subsection~\ref{artificiel}, we illustrate this difference with an artificial but interesting example, 
for which the reduced LAP \eqref{eq:tal-reduit} holds true and the usual Mourre estimate \eqref{eq:mourre2} is false. This example is however related to the situation in \cite{dj}. 

Paper's organisation: In Section~\ref{nouv-app}, we introduce the main tools of our new 
approach. Admitting Theorem~\ref{th:main}, we prove Theorem~\ref{th:mourre} in Subsection~\ref{local}. 
Section~\ref{new-proof} is devoted to the proof of Theorem~\ref{th:main}. In Section~\ref{etude-tal-reduit}, 
we prove Theorems~\ref{th:mourre2} and~\ref{th:mourre2'} on the reduced resolvent. Technical tools 
are collected in Appendices~\ref{calcul-symb}, ~\ref{dev-commut}, and~\ref{esti-techn}. 

{\bf Acknowledgments:} We would like to thank Jan Derezi\'nski, Vladimir Georgescu, Jacob S. 
M\o ller, and Francis Nier for helpful discussions. The two authors were partially supported by the contract MERG-CT-2004-006375 funded 
by the European commission and the first one was also supported by the Postdoctoral Training Program
HPRN-CT-2002-0277. The second author thanks the members of the Institut of Mathematics ``Simion Stoilow" of the 
Romanian Academy in Bucharest for their kind hospitality.

\section{A new approach of the LAP.} \label{nouv-app}

We explain in this section our strategy to prove Theorem~\ref{th:main}
below, a stronger version of Theorem~\ref{th:mourre}.

\subsection{Basic facts and notation}\label{base}

In this subsection, we introduce some notation and recall known basic results. We refer to 
\cite{abg} for details.

In the text, we use the letter $\Ic$ to denote an interval of
$\R$. For such a $\Ic$, we denote by $\overline{\Ic}$
(resp. $\mathring{\Ic}$) its closure (resp. its interior). 
The scalar product $\langle \cdot ,\cdot \rangle$ in $\Hr$ is right linear and $\|\cdot\|$ denotes 
the corresponding norm and also the norm of bounded operators on $\Hr$. 
If $T$ is a bounded operator on $\Hr$ and $k\in \N$, we say that $T\in \Cc^k(A)$ 
if, for all $f\in \Hr$, the map $\R\ni t\mapsto e^{itA}Te^{-itA}f\in \Hr$ has the usual 
$\Cc^k$ regularity. It turns out that $T\in \Cc^k(A)$ if and
only if, for a $z$ outside the spectrum of $T$, $(T-z)^{-1}\in \Cc^k(A)$. For such $T$,
$T\in \Cc^1(A)$ if and only if the form $[T,A]$ defined on $\Dc(A)\times \Dc(A)$ extends to 
a bounded operator $\ad_A^1(T)=[T,A]$ if and only if $T$ preserves $\Dc (A)$ and 
the operator $TA-AT$, defined on $\Dc(A)$, extends to a bounded operator in $\Hr$. Furthermore 
$T\in \Cc^k(A)$ if and only if the iterated commutator $\ad_A^p(T):=[\ad_A^{p-1}(T),A]$ 
are bounded for $p\leq k$. In particular, for $T\in \Cc^1(A)$, $T\in \Cc^2(A)$ if
and only if $[T,A]\in \Cc^1(A)$. For unbounded self-adjoint operator, we defined the
$\Cc^k(A)$ regularity in Section~\ref{intro}. Let $H$ is (unbounded) self-adjoint operator and 
$\Ic$ a bounded interval. Recall that $E_\Ic (H)$ denotes the spectral projection of $H$ above $\Ic$. 
If $H\in \Cc^1(A)$ then the form $[H,iA]$ defined on $(\Dc(H)\cap\Dc(A))\times (\Dc(H)\cap\Dc(A))$ extend to a bounded operator from $\Dc(H)$ to its dual for the graph norm. In particular, \eqref{eq:mourre2} makes sense. 
A justification of Remark~\ref{r:K=0} can be found in \cite{abg} but 
we give it in the proof of Theorem~\ref{th:mourre2} (for $P=0$). 
The following propositions and remark will be useful later. 
\begin{proposition}\label{p:rang1} 
For $f,g\in \Dc (A)$, the finite rank operator $|f\rangle \langle
g|:h\rightarrow \langle g,h\rangle \cdot f$ belongs to $\Cc^1(A)$ and 
$[|f\rangle \langle g|,A]=|f\rangle \langle Ag|-|Af\rangle \langle
g|$. In particular, if $f,g\in \Dc (A^2)$, $|f\rangle \langle
g|\in \Cc^2(A)$. If $P$ is a finite rank projection, the range of which is included in 
$\Dc (A^k)$ with $k\in \N$, then $P\in \Cc ^k(A)$
\end{proposition}  
\proof Since $R:=|f\rangle \langle g|$ preserves $\Dc(A)$, the commutator $[|f\rangle \langle g|,A]$ 
is well defined on $\Dc(A)$ and $[|f\rangle \langle g|,A]=|f\rangle \langle Ag|-|Af\rangle \langle g|$, which extends to a bounded operator. Thus $R\in\Cc^1(A)$. Applying the first result to $[R,A]$, $|f\rangle \langle g|\in \Cc^2(A)$ if $f,g\in \Dc (A^2)$. Since  $P=\sum _{1\leq n\leq N}|
f_n\rangle \langle f_n|$ with $N\in\N$, an induction argument gives the last result. \qed
\begin{proposition}\label{p:fermeture} 
Let $(T_n)_n$ be a sequence of bounded operators such that, $T_n\in \Cc^1(A)$, for all $n$, 
and such that there exist bounded $S,T$ such that $T_n\to T$ and $[T_n,A]\to S$ in the norm 
topology. Then $T\in \Cc^1(A)$ and $S=[T,A]$. 
\end{proposition}  
\proof See Lemma 2.5 in \cite{ggm1}. \qed
\begin{rem}\label{r:TAL-s}
The LAP, respectively to $(\Ic,0,A)$, holds true for $H$ if and only if $H$ has no spectrum in 
$\Ic$. The LAP for $H$, respectively to $(\Ic,s,A)$, implies the LAP for $H$, respectively to 
$(\Ic,s',A)$, for any $s'\geq s$. For $H=-\Delta$ the Laplace operator in $\R^d$ and 
$A$ the multiplication operator by $\langle x\rangle$, it is known 
that LAP for $H$, respectively to $(\Ic,s,A)$, holds true if and only if $s>1/2$ (cf. \cite{h}). 
\end{rem}

\subsection{Local regularity and main result} \label{local}
In Theorem~\ref{th:mourre}, the LAP \eqref{eq:lap} and the 
Mourre estimate \eqref{eq:mourre2} are localized in $H$. 
It is quite natural to try to replace $H$ and the global assumption 
$H\in \Cc^2(A)$ by some local version. By \cite{abg}, we have 
\begin{proposition}\label{p:reguv} 
Let $\varphi\in\Cc_c^\infty(\R)$. Suppose $H\in\Cc^k(A)$ for a certain 
$k\in\N$. Then, $\varphi(H)\in \Cc^k(A)$.  
\end{proposition}  

For any $\tau\in\Cc^\infty_c(\R)$, we define the bounded operator 
\begin{equation}\label{eq:def-H-tau}  
H_\tau := H\tau (H). 
\end{equation}  
It turns out that we can deduce the LAP for $H$ respectively to $(\Ic,s,A)$ 
from the LAP for $H_\tau$ respectively to $(\Ic,s,A)$, if $\tau =1$ near 
$\Ic$, as seen in Proposition~\ref{p:equi} below. Thus $H_\tau$ is a good local 
(and bounded) version of $H$. From \cite[Proposition 2.1]{s}, we pick the following  
\begin{lemma}\label{l:mourrer} 
Let $\Ic$ be bounded, open interval. Suppose that $H\in \Cc^1(A)$ and that the Mourre 
estimate \eqref{eq:mourre2} holds true on $\Ic$. Take $\theta\in\Cc^\infty_c(\Ic)$ and 
$\tau\in\Cc^\infty_c(\R)$ such that $\tau\theta=\theta$. Then $H_\tau\in \Cc^1(A)$ and  
\begin{eqnarray}\label{eq:mourrer} 
\theta(H)[H_\tau, iA]\theta(H)\geq c \theta^2(H)+\theta(H)K\theta(H). 
\end{eqnarray}  
\end{lemma}  
\proof By Proposition~\ref{p:reguv}, $H_\tau\in \Cc^1(A)$. 
For $f\in\Dc(A\theta(H))$, 
\begin{eqnarray*} 
\langle H\theta(H)f, iA\theta(H)f\rangle -   \langle iA \theta(H)f, 
H\theta(H)f \rangle\geq c \|\theta(H)f\|^2 + \langle f, Kf \rangle. 
\end{eqnarray*}  
Now, use that $H\theta(H)= H\tau(H)\theta(H)$. Finally, $\Dc(A\theta(H))$ is 
dense in $\Hr$ since $\theta(H)A$ is closed with a dense domain. \qed 
 
\begin{rem} In general, one should not expect a ``real'' Mourre estimate for 
$H_\tau$ of the form  
\begin{eqnarray*} 
\varphi (H_\tau)[H_\tau, iA]  \varphi (H_\tau)\geq c \varphi ^2(H_\tau)+K, 
\end{eqnarray*} 
for a certain function $\varphi $ which satifies the same hypothesis as 
$\theta$ in Lemma~\ref{l:mourrer}. Indeed, since $0\in \supp \theta$, there is no such function 
$\varphi $ such that $\varphi (t\tau(t))=\theta(t)$ for all $t\in\R$. 
\end{rem} 
Given an open interval $\Ic$ and $k\in\N$, we say that $H$ is {\em locally
  of class $\Cc^k(A)$ on $\Ic$}, we write $H\in \Cc^k_\Ic(A)$, if, for 
all $\varphi\in\Cc_c^\infty(\Ic)$, 
$\varphi(H)\in \Cc^k(A)$. This is a local version of the regularity $\Cc^k(A)$ which 
was already used in \cite{s}. 

{\bf Proof of Theorem~\ref{th:mourre}:} Let $\Ic ''$ be open such that 
$\overline{\Ic}\subset \Ic ''$. By Lemma~\ref{p:reguv}, $H\in \Cc^2_{\Ic ''}(A)$. Let $\tau \in
\Cc^\infty_c(\Ic '')$ such that $\tau =1$ near $\Ic$. Let $\Ic _1$ be closed such that 
$\Ic '\subset \mathring{\Ic _1}$ and $\Ic_1\subset \Ic$. 
Let $\theta\in\Cc^\infty_c(\Ic)$ such that $\theta =1$ on $\Ic _1$. By Lemma~\ref{l:mourrer} 
and \eqref{eq:mourre2}, we derive \eqref{eq:mourrer}, which implies 
\begin{eqnarray*} 
E_{\Ic _1}(H)[H_\tau, iA]E_{\Ic _1}(H)\geq c E_{\Ic _1}(H)+0, 
\end{eqnarray*} 
since $\theta =1$ on $\Ic _1$. 
Thus Theorem~\ref{th:main} below applies yielding the LAP for $H$ respectively to 
$(\Ic ',s,A)$.  \qed

So the proof of Theorem~\ref{th:mourre} reduces to the proof of the following 
stronger result, which is our main result. 
 
\begin{theorem}\label{th:main} 
Let $\Ic$ be a bounded, open interval. Let $\Ic ''$ be an open interval 
such that $\overline{\Ic}\subset \Ic ''$. 
Let $H\in \Cc^2_{\Ic ''}(A)$ and $\tau\in\Cc^\infty_c(\Ic '')$ 
such that $\tau=1$ near $\Ic$. Suppose the strict Mourre estimate 
\begin{eqnarray}\label{eq:mourremain} 
E_\Ic (H)[H_\tau, iA]  E_\Ic(H)\geq c E_\Ic(H),\mbox { with }  c\g 0,
\end{eqnarray}  
holds  true. Then, for any $s\g 1/2$ and any compact interval $\Ic'$
with $\Ic'\subset \mathring{\Ic}$, 
the LAP respectively to $(\Ic', s, A)$ holds true for $H_\tau$ and $H$. 
\end{theorem} 
\proof See Subsection~\ref{s:proof-main} (and Subsection~\ref{s:sketch} for a sketch). 
\qed

\begin{rem}\label{r:s}
In \cite{s}, the previous result is proved under a weaker local
regularity assumption (slightly stronger than $\Cc^1_{\Ic ''}(A)$),
using Mourre's differential inequality technics. Furthermore, an
example of multiplication operator $H$ and of conjugate operator $A$
is given such that $H\not\in \Cc^1(A)$ but $H\in \Cc^1_{\Ic}(A)$, 
for some $\Ic$. 
\end{rem}

\begin{rem}\label{r:exemple}
Assume that Theorem~\ref{th:mourre} applies to some operators $H$ and $A$ on some interval $\Ic$. 
Let $\Ic ''$ be open such that $\overline{\Ic}\subset \Ic ''$. Let $\varphi:\R\longrightarrow\R$ 
be a borelian, increasing function such that, for all $t\in \Ic''$, $\varphi (t)=t$. Then 
Theorem~\ref{th:main} applies with $H$ replaced by $\varphi (H)$. Since $\varphi$ may be irregular outside 
$\Ic ''$, we do not know if $\varphi (H)\in \Cc^1(A)$, so if Theorem~\ref{th:mourre} applies to $\varphi (H)$. 
\end{rem}

\subsection{Special sequences and the LAP} 
In this subsection, we introduce our main tool and its properties. We proceed 
like in \cite{jec} and use the terminology appearing in this semi-classical 
setting.  
 
\begin{definition} \label{suite-speciale}
A \emph{special sequence} $(f_n, z_n)_n$ for $H$ associated to $(\Ic,s, A)$, as in \eqref{eq:lap}, 
is a sequence $(f_n,z_n)_n\in (\Dc(H)\times \C)^\N$ such that, for certain $\lambda\in \Ic$ 
and $\eta \geq 0$,  $\Ic\ni\re(z_n)\rightarrow\lambda$, $0\neq\im(z_n)\rightarrow 0$, 
$\|\langle A\rangle ^{-s} f_n \|\rightarrow \eta$, $(H-z_n)f_n\in\Dc(\langle A\rangle^{s})$, 
and $\|\langle A\rangle ^{s}(H-z_n)f_n \|\rightarrow 0$. The limit $\eta$ is called the 
\emph{mass} of the special sequence. 
\end{definition}  
We give the link between this notion and the LAP in 
\begin{proposition}\label{p:non-lap}  
Given $s\geq 0$ and a compact interval $\Ic$, the LAP for $H$ respectively to $(\Ic, s, A)$ is false 
if and only if there exists a special sequence $(f_n, z_n)_n$ for $H$ associated to $(\Ic, s, A)$ with a 
positive mass.  
\end{proposition}  
  
\proof   
Suppose the LAP to be false. There exist a sequence $(k_n)_n$ of nonnegative numbers, going to infinity, 
a sequence $(g_n)_n$ of non-zero elements of $\Hr$, and a sequence $(z_n)_n$ of complex numbers 
such that $\re (z_n)\in \Ic$, $0\neq \im (z_n)\to 0$, and    
\begin{equation}\label{eq:taille-resolv}  
\bigl\|\langle A\rangle  ^{-s}(H-z_n)^{-1}\langle A\rangle  ^{-s}g_n\bigr\|\ = 
\ k_n\, \|g_n\|\ = 1. 
\end{equation}  
Setting $f_n=(H-z_n)^{-1}\langle A\rangle  ^{-s}g_n$, $f_n\in \Dc(H)$, 
$(H-z_n)f_n\in \Dc(\langle A\rangle^{s}) $, and, by 
(\ref{eq:taille-resolv}),   
\begin{eqnarray*}
\bigl\|\langle A\rangle  ^{-s}f_n\bigr\|\ = \ 1 \ \mbox{and}\  \,   
\|\langle A\rangle ^{s}(H-z_n)f_n\|\ = 1/k_n\to 0.  
\end{eqnarray*} 
Up to a subsequence, we can assume that $\re (z_n)\to \lambda\in \Ic$. Now, we assume the 
LAP true and consider $(f_n, z_n)_n$, a special sequence for $H$ associated to
$(\Ic, s, A)$. By \eqref{eq:lap}, there exists $c>0$ such that 
\begin{eqnarray*}
\bigl\|\langle A\rangle  ^{-s}f_n\bigr\|\ \leq \ c \|\langle A\rangle ^{s}(H-z_n)f_n\|.  
\end{eqnarray*}
This implies $\eta =0$. \qed 
 
The previous result can be partially localized in energy.  

\begin{proposition}\label{p:easy} 
Let $(\Ic, s, A)$ be a triplet as in \eqref{eq:lap} with $0\leq s\l 1$. Let $\Ic ''$ be open such 
that $\overline{\Ic}\subset \Ic ''$ and $H\in\Cc^1_{\Ic ''}(A)$. Let $\theta\in\Cc^\infty_c(\R)$ 
such that $\theta =1$ near $\Ic$. Let $\varphi :\R\rightarrow\R$ borelian such that,  
for $t\in \supp\theta$, $\varphi(t)=t$. Let $(f_n, z_n)_n$ be a special sequence for $H$ associated to 
$(\Ic, s, A)$ with mass $\eta$. Then, writing $\tilde\theta =1-\theta$,  
\begin{enumerate} 
\item $\tilde\theta(H) f_n$ tends to $0$, 
\item $(\theta(H)f_n, z_n)_n$ is a special sequence for $\varphi(H)$ associated to $(\Ic, s, A)$ 
with mass $\eta$. 
\end{enumerate} 
\end{proposition} 
 
\proof Since 
$$\|\tilde{\theta}(H)f_n\|\ \leq \ \|\tilde{\theta}(H)(H-z_n)^{-1} 
\langle A \rangle^{-s}\|\, \cdot \, \| \langle A
\rangle^{s}(H-z_n)f_n\| 
$$ 
and since $t\mapsto \tilde{\theta}(t)/(t-z_n)$ is uniformly bounded in
$n$, $\|\tilde{\theta}(H)f_n\|$ tends to $0$.  
Since $s\geq 0$, $\|\langle A \rangle^{-s}\tilde{\theta}(H)f_n\|\rightarrow 0$ and therefore  
$\|\langle A \rangle^{-s}\theta(H)f_n\|\rightarrow \eta$. Since $H\in\Cc^1_{\Ic ''}(A)$, 
$\theta(H)\in \Cc^1(A)$. Since $s\l 1$, $\|\langle A \rangle^{s} 
\theta(H) \langle A \rangle^{-s}\|$ is bounded, by Proposition~\ref{p:regu}. Now, 
$$\|\langle A \rangle^{s} (\varphi(H)-z_n) \theta(H)f_n \|\leq  
\|\langle A \rangle^{s} \theta(H) \langle A \rangle^{-s}\|\cdot 
\|\langle A \rangle^{s} (H-z_n)f_n \|$$ 
which tends to $0$.  \qed 

Now we can perform the reduction to some $H_\tau$ (cf. \eqref{eq:def-H-tau}). 
  
\begin{proposition}\label{p:equi} 
Let $(\Ic, s, A)$ be a triplet as in \eqref{eq:lap} with $0\leq s\l 1$. Let $\Ic ''$ be open such that 
$\overline{\Ic}\subset \Ic ''$ and $H\in\Cc^1_{\Ic ''}(A)$. Let $\tau\in\Cc^\infty_c(\Ic '')$ such that 
$\tau =1$ near $\Ic$. If the LAP respectively to $(\Ic,s,A)$ holds true for $H_\tau$ then it holds
true for $H$.  
\end{proposition} 
\proof By contraposition, the result follows from Propositions~\ref{p:non-lap} and~\ref{p:easy}.  \qed 

\begin{rem} There is another proof of Proposition~\ref{p:equi}. Let $\theta\in \Cc_c^\infty(\R)$ 
with $\theta =1$ near $\Ic$ and $\tau\theta=\theta$. Then, using a Neumann serie 
for $|z|$ large enough and $z\not\in \R$, we can show that $(H-z)^{-1}\theta (H)=
(H_\tau-z)^{-1}\theta (H)$. By analyticity, this holds true for $z\not\in \R$. 
Therefore, if the LAP respectively to $(\Ic,s,A)$ is true for $H_\tau$ 
so is it for $H$, since $\langle A \rangle^{s} 
\theta(H) \langle A \rangle^{-s}$ is bounded. 
\end{rem}

\subsection{A Virial-like Theorem.}

In Remark~\ref{r:viriel+mourre}, we recalled the Virial Theorem. 
Our approach is based on the following Virial-like result. 

\begin{proposition}\label{p:virial} 
Let $(f_n,z_n)_n$ be a special sequence for a \emph{bounded}, self-adjoint operator $H_b$ 
respectively to $(\Ic, s, A)$, as in \eqref{eq:lap} with $s\geq 0$. For any \emph{bounded} borelian function $\phi$, 
\begin{eqnarray*} 
\lim_{n \rightarrow \infty} \langle f_n, [H_b, \phi(A)] f_n\rangle =0. 
\end{eqnarray*}  
\end{proposition}  
\proof  
Since $[H_b, \phi(A)]=[H_b-z_n, \phi(A)]$, 
\begin{eqnarray*} 
\langle f_n, [H_b, \phi(A)] f_n\rangle
&=&2i\im(z_n) \langle f_n,\phi(A) f_n\rangle \\ 
&& \hspace*{-2cm} + \langle (H_b-z_n)f_n, \phi(A) f_n \rangle + 
\langle \phi(A)^* f_n, (H_b-z_n) f_n \rangle. 
\end{eqnarray*} 
By Definition~\ref{suite-speciale}, there exists $C>0$ such that 
\begin{eqnarray*} 
|\langle (H_b-z_n) f_n, \phi(A) f_n \rangle| &\leq& 
|\langle \langle A\rangle^{s}  (H_b-z_n) f_n, \langle A\rangle^{-s} 
\phi(A) f_n \rangle | \\ 
&\leq& C \|\phi (A)\|\cdot \|\langle A\rangle^{s}  (H_b-z_n) f_n\| \build{\rightarrow}_{n 
\rightarrow \infty}^{} 0. 
\end{eqnarray*}  
Similarly, $\lim\langle \phi(A)^\ast f_n ,(H_b-z_n) f_n\rangle=0$.  
By Definition~\ref{suite-speciale},  
\begin{eqnarray*} 
 \im(z_n)\cdot\| f_n\|^2&=& \im \langle f_n, (H_b- z_n) f_n \rangle \\ 
&=& \im \langle \langle A \rangle^{-s} f_n, \langle A \rangle^{s} (H_b- z_n) f_n  
\rangle \build{\rightarrow}_{n \rightarrow \infty}^{} 0.  
\end{eqnarray*}  
Since 
\begin{eqnarray*} 
|\im(z_n) \langle f_n,\phi(A) f_n\rangle|\leq |\im(z_n)|\cdot\| f_n\|^2 \cdot\|\phi(A)\|, 
\end{eqnarray*}  
we obtain the desired result. \qed 

\begin{rem} If $H_b$ is not bounded, Proposition~\ref{p:virial} works, provided the commutator 
$[H_b, \phi(A)]$ is considered as quadratic form. 
\end{rem}

\subsection{Sketch of our proof and interpretation.} \label{s:sketch}

To prove Theorem~\ref{th:main}, we only need to show the LAP for $H_\tau$ on $\Ic'$ by 
Proposition~\ref{p:equi}. In view of Proposition~\ref{p:non-lap}, we consider 
a special sequence $(f_n, z_n)_n$ for $H_\tau$ associated to the triplet $(\Ic, s, A)$, 
with $s\g 1/2$, and we show that $\eta=0$. By Remark~\ref{r:TAL-s}, we
may assume that $s\in ]1/2;2/3[$. For $R>1$, let $\cchi _R+\tilde{\cchi}_R=1$ be 
a smooth partition of unity on $\R$ with $\cchi _R$ localized in $\{t\in \R;|t|\leq 2R\}$. 
It suffices to show that $\lim_{R\to\infty}\limsup_{n\rightarrow\infty}\|\tilde{\cchi}_R(A)
\langle A\rangle ^{-s}f_n\|=0$ and $\lim_{R\to\infty}\limsup_{n\rightarrow\infty}\|
\cchi _R(A)\langle A\rangle^{-s}f_n\|=0$. From the strict Mourre estimate 
\eqref{eq:mourremain}, we deduce \eqref{eq:mourrer} with $K=0$. 
We apply the latter to $\tilde{\cchi}_R(A)\langle A\rangle ^{-s}f_n$. After several commutations, 
the use of Proposition~\ref{p:virial}, and the use of the assumption $s\g 1/2$, we find some 
$\epsilon>0$ such that, for all $R>1$,  
\begin{eqnarray}\label{eq:controle-infini} 
\limsup_{n\rightarrow\infty}\|\tilde{\cchi}_R(A)\langle A\rangle ^{-s}f_n\|=O(R^{-\epsilon}). 
\end{eqnarray} 
Next we apply the Mourre estimate \eqref{eq:mourrer} to $\cchi _R(A)f_n$. After several commutations, 
the use of Proposition~\ref{p:virial}, and the use of \eqref{eq:controle-infini}, we get
$\lim_{R\to\infty}\limsup_{n\rightarrow\infty}\|\cchi _R(A)f_n\|=0$. Since $s\geq 0$, we obtain 
the desired results yielding $\eta=0$. 

This proof provides the following new interpretation of Theorems~\ref{th:mourre} and~\ref{th:main}. 
The strict Mourre estimate excludes the existence of a special sequence of positive mass, yielding the LAP, 
in a similar way as it excludes the existence of a bound state in Remark~\ref{r:viriel+mourre}. 
Our Virial-like Theorem plays the role of the usual Virial Theorem.

\section{A new proof of the LAP.}\label{new-proof}

Here we complete the proof of Theorem~\ref{th:main} sketched in Subsection~\ref{s:sketch}. 
We assume the assumptions of Theorem~\ref{th:main} satisfied and take
some interval $\Ic'\subset \mathring{\Ic}$. Let $\theta\in\Cc_c^\infty(\Ic)$ with $\theta =1$ on $\Ic'$. 
Applying $\theta (H)$ on both sides, we deduce from \eqref{eq:mourremain} the strict Mourre 
estimate \eqref{eq:mourrer} (i.e. with $K=0$). We consider a special sequence $(f_n, z_n)_n$ 
for $H_\tau$ associated to $(\Ic, s, A)$ with $s\in ]1/2;2/3[$. By Proposition~\ref{p:easy}, we 
may assume that $\theta (H)f_n=f_n$, for all $n$. Let us fix some notation. Let 
$\cchi\in \Cc^\infty_c(\R)$ such that  
\begin{eqnarray}\label{eq:chi} 
\cchi=1 \mbox{ on } [-1,1] \mbox{ and } \cchi=0 \mbox{ on } \R\setminus [-2,2]. 
\end{eqnarray} 
We shall require other properties satisfied by $\cchi$ (see \eqref{eq:chi2} below). 
For $R\g 1$, we set $\cchi_R(x)=\cchi(x/R)$ for all $x\in \R$ and 
$\tilde\cchi_R=1-\cchi_R$.  We denote by $O_R(\cdot)$ (resp. $o_R(\cdot)$) the Landau symbol 
$O$ (resp. $o$) where the subscript $R$ means that the bound (resp. the limit) is uniform w.r.t. 
the other variables.

\subsection{Proof of Theorem \ref{th:main}} \label{s:proof-main}
Let $\cchi\in \Cc^\infty_c(\R)$ satisfying \eqref{eq:chi} and \eqref{eq:chi2}.  
From Proposition \ref{p:mass2} and Corollary \ref{p:mass} below, we derive that, for all 
$\varepsilon \g 0$, 
\begin{eqnarray*} 
\eta &=& \lim_{n\rightarrow \infty}\|\langle A \rangle^{-s} f_n\| \\ 
&\leq&  \limsup_{n\rightarrow \infty}\, (\|\tilde\cchi_R(A) \langle A 
\rangle^{-s} f_n \| +  \|\cchi_R(A) \langle A \rangle^{-s} f_n \|) \\
&\leq&  \limsup_{n\rightarrow \infty}\, (\|\tilde\cchi_R(A) \langle A 
\rangle^{-s} f_n \| +  \|\cchi_R(A)f_n \|)= O(R^{2s-2+\varepsilon}). 
\end{eqnarray*}  
Letting $R$ go to infinity, we obtain that $\eta =0$. By Proposition 
\ref{p:non-lap}, the LAP holds true for $H_\tau$ respectively to $(\Ic', s, 
A)$. \qed

\subsection{A ``large $|A|$" estimate. }  
We stress that, in this subsection, we suppose that $1/2\l s\l 1$.  
The aim of this part is to show 
\begin{proposition}\label{p:massest} 
Let $\Ic'$ be closed with $\Ic'\subset\mathring{\Ic}$ and let $(f_n,z_n)_n$ be a special sequence for $H_\tau$ 
respectively to $(\Ic', s, A)$ with $1\g s\g 1/2$. Assume that the 
Mourre estimate \eqref{eq:mourrer} holds true with $\theta=1$ on $\Ic'$ and $K=0$.  
Let $\cchi\in \Cc^\infty_c(\R)$ satisfying \eqref{eq:chi} and \eqref{eq:chi2} (below). 
Then, there exist $c'>0$, $R_1>2$, and a family $(\phi_R)_{R>1}$ in $L^\infty(\R)$, 
such that, for all $R\geq R_1$, 
\begin{eqnarray}\nonumber 
\langle f_n, [H_\tau, i\phi_R(A)] f_n\rangle&\geq& c' \|\tilde\cchi_R(A) 
 \langle A \rangle^{-s} f_n \|^{2} \\ \label{eq:massest} 
&\hspace*{-2cm}+&\hspace*{-1cm}O_R(R^{-1})\cdot  
 \|\tilde\cchi_{R}(A)\langle A  \rangle^{-s} f_n  \| \\ \nonumber 
&\hspace*{-2cm}+&\hspace*{-1cm}O_R(R^{2s-2})\cdot 
 \|\tilde\cchi_{R/2}(A)\langle A  \rangle^{-s} f_n  \| .     
\end{eqnarray} 
\end{proposition}   
 
\begin{corollary}\label{p:mass} 
Under the hypotheses of Proposition \ref{p:massest}, 
\begin{eqnarray}\label{eq:mass} 
\forall \alpha >2s-2\, ,\hspace{.5cm}\build{\limsup}_{n \rightarrow \infty}^{}\|\tilde\cchi_R(A) \langle A 
\rangle^{-s} f_n \|= O(R^{\alpha}). 
\end{eqnarray} 
\end{corollary}  
\proof Note first that, for $a>0$, $\varepsilon\geq 0$, and $b,c\in \R$, 
\begin{eqnarray}\label{eq:carre} 
\quad \quad \quad\varepsilon\geq aX^2 +b X + c^2\implies  
|X|\leq \sqrt{\varepsilon /a }+ O(|b|+|c|),    
\end{eqnarray} 
the latter term being independent of $\varepsilon$. Since $(\|\langle A 
\rangle^{-s} f_n \|)_n$ is bounded by Definition~\ref{suite-speciale}, it suffices to 
prove \eqref{eq:mass} for large $R$. For fixed $R\geq R_1$, we combine 
\eqref{eq:massest} with \eqref{eq:carre} and Proposition~\ref{p:virial} to get 
\begin{eqnarray}
&&\limsup_{n\rightarrow \infty}\|\tilde\cchi_R(A)\langle A \rangle^{-s}
f_n \| \leq O(R^{-1}) \nonumber   \\
&&+ \, O(R^{s-1}) 
\cdot \limsup_{n\rightarrow \infty}\|\tilde\cchi_{R/2}(A)\langle A
\rangle^{-s} f_n  \|^{1/2}.  \label{eq:limsup}
\end{eqnarray} 
We use a bootstrapping argument. Since $(\|\langle A\rangle^{-s} f_n \|)_n$ is bounded, so 
is $(\|\tilde\cchi_{R/2}(A)\langle A\rangle^{-s} f_n\|)_n$. Then \eqref{eq:limsup} gives 
\eqref{eq:mass} for $\alpha=\alpha_0=s-1$.  Now we use this new estimate in \eqref{eq:limsup} 
to get \eqref{eq:mass} for $\alpha=\alpha_1=3(s-1)/2$. By induction, we get \eqref{eq:mass} for 
a sequence $(\alpha_n)_n$ satisfying $\alpha_{n+1}= \alpha_n/2 +(s-1)$, for all $n\in \N$. 
By a fixed point argument, $\alpha _n\to 2(s-1)$. This yields the result.  \qed 

Our strategy to prove Proposition \ref{p:massest} is the following. We apply the strict Mourre 
estimate \eqref{eq:mourrer} (with $K=0$) to $\tilde\cchi_R(A)\langle A\rangle^{-s} f_n$. 
We move the $\theta (H)$ to the $f_n$, which absorb them, since $\theta (H)f_n=f_n$. We want 
to pull the weights $\tilde\cchi_R (A)\langle A\rangle^{-s}$ into the commutator $[H_\tau,A]$, 
in order to get the term on the l.h.s of \eqref{eq:massest} with $\phi _R(t)=
t\langle t\rangle ^{-2s}\tilde\cchi_R(t)^2$. In view of the proof of Corollary~\ref{p:mass}, we 
need $2s\geq 1$. Our manipulation produces of course error terms which should be small. Using 
$s>1/2$, we actually prove this smallness if we only move the $\theta (H)$ and the 
$\langle A\rangle^{-s}$. Choosing appropriate functions $\phi _R$, we can move the 
$\tilde\cchi_R(A)$ into the commutator producing an error term which has the good sign, 
up to small enough terms. To this end, we choose more carefully the function $\cchi$ in 
\eqref{eq:chi}. We demand that $\cchi$ satisfies \eqref{eq:chi} and that  
\begin{eqnarray}\label{eq:chi2} 
\tilde\cchi:=1-\cchi=\tilde\cchi_+ +\tilde\cchi_-,  
\end{eqnarray} 
where $\tilde\cchi\ind_{\R^\sigma}= \tilde\cchi_\sigma$ for $\sigma\in\{-,+\}$, such 
that $\tilde\cchi_\sigma$ and $\sigma\tilde\cchi_\sigma'$ are square of some smooth functions 
(see for instance the appendix in \cite{dg} for their existence). Let $R>1$. We set 
$\tilde\cchi_{\sigma, R}=\tilde\cchi_\sigma(t/R)$. Notice that $\tilde\cchi_R= 
\tilde\cchi_{+, R}+\tilde\cchi_{-, R}$ and $\tilde\cchi_R^2= \tilde\cchi_{+, R}^2
+\tilde\cchi_{-, R}^2$. 
\proof[Proof of Proposition \ref{p:massest}] Let 
\begin{eqnarray}\label{eq:h} 
h&\g &\sup_{t\in\R} |t|\langle t\rangle^{-2s}.
\end{eqnarray}  
Let $R>1$ and $\phi _R\in \Cc^\infty(\R)$ defined by  
\begin{eqnarray}\label{eq:phi-R} 
\phi_R(t)& =& \sum_{\sigma\in\{-,+\}}\tilde\cchi_{\sigma, R}^2(t)\bigl(\sigma h - t 
\langle t \rangle^{-2s}\bigr).    
\end{eqnarray}  
For all $f\in \Hr$, 
\begin{eqnarray}\nonumber 
\langle f, [H_\tau, i \phi_R(A)] f\rangle&=& 
\langle f, \tilde\cchi_R(A)[H_\tau, - i A\langle A\rangle^{-2s}] 
\tilde\cchi_R(A)f \rangle \\\label{eq:step} 
   &&\hspace*{-2cm}+\,\,\, 2\!\sum_{\sigma\in\{-,+\}}\re \langle f, (\sigma h 
- A \langle A\rangle^{-2s}) \tilde\cchi_{\sigma, R}(A) [H_\tau, i 
  \tilde\cchi_{\sigma, R}(A)] f \rangle.  
\end{eqnarray}  
We can find some $R_1>2$ (see Lemma~\ref{l:mino-commut} below)  such that, for $R\geq R_1$,   
\begin{eqnarray}\nonumber 
\langle \tilde\cchi_R(A) f_n, [H_\tau,-iA\langle A\rangle^{-2s}] 
\tilde\cchi_R(A) f_n\rangle\geq c' \|\tilde\cchi_R(A)\langle A\rangle^{-s} 
f_n\|^2\\\label{eq:1}  
&\hspace*{-12cm}&\hspace*{-6cm}+O_R(R^{-1})\cdot 
\|\tilde\cchi_R(A)\langle A\rangle^{-s} f_n\| 
\end{eqnarray} 
with $c'= (2s-1)^{-1}c/2>0$. 
Since we are not able to estimate properly the second term on the r.h.s of \eqref{eq:step}, we indend 
to use  some positivity argument to get rid of it. In view of \eqref{eq:egalite}, we choosed $\phi_R$ 
in \eqref{eq:phi-R} such that, the function $\psi _R\in \Cc_c^\infty(\R)$ defined by 
\begin{eqnarray}
\psi_R(t)&=&R\bigl(\phi_R'(t)-\tilde\cchi_R^2(t)(d/dt)(-t\langle
t\rangle^{-2s})\bigr) \nonumber\\
&=& \sum_{\sigma\in \{-,+\}}(\tilde\cchi_{\sigma}')_R(t) (\sigma h- t\langle t 
\rangle^{-2s})\tilde\cchi_{\sigma, R}(t)\label{eq:positive} 
\end{eqnarray}  
is the square of a smooth function. We put a factor $R$ in front to ensure that the family 
$(\psi_R)_R$ is bounded in some symbol space (see Lemma~\ref{l:step}). 
Notice that $\supp \psi_R\subset [-2R, -R]\cup [R,2R]$. 
We define  
\begin{eqnarray}\label{eq:cr} 
C_R:=\psi_R^{1/2}(A)\hspace{.2cm}\mbox{and note that}\hspace{.2cm}C_R\tilde\cchi_{R/2}(A)= C_R.  
\end{eqnarray}  
We can show (see Lemma~\ref{l:mino-C-R} below) that 
\begin{eqnarray} \label{eq:mino-C-R}
\hspace{.6cm}\langle C_R f_n, [H_\tau, iA] C_R f_n \rangle  &\geq &O_R(R^{2s-1}) 
\cdot \|\tilde\cchi_{R/2}(A) \langle A\rangle^{-s} f_n\|. 
\end{eqnarray} 
By Lemma \ref{l:C}, 
\begin{eqnarray*} 
\langle C_R^2 f_n, [H_\tau, iA]  f_n \rangle  &\geq& O_R(R^{2s-1})\cdot \|\tilde\cchi_{R/2}(A) 
\langle A\rangle^{-s} f_n\|,  
\end{eqnarray*} 
since $\|\tilde\cchi_{R/2}(A) \langle A\rangle^{-s} f_n\|=O_R(R^0)$. 
Now, by Lemma \ref{l:exp}, 
\begin{eqnarray*} 
\langle f_n, (\sigma h - A 
\langle A\rangle^{-2s}) \tilde\cchi_{\sigma, R}(A)[H_\tau , \tilde\cchi_{\sigma, 
    R}(A)] f_n \rangle&\geq&\\ 
&\hspace*{-10cm}+&\hspace*{-5cm}O_R(R^{2s-2})\cdot\|\tilde\cchi_{R/2}(A)\langle A\rangle^{-s}f_n\|. 
\end{eqnarray*}  
This yields, together with \eqref{eq:step} and \eqref{eq:1}, the
result. \qed  

To prove \eqref {eq:1} and \eqref {eq:mino-C-R}, we need the following lemmata. 
\begin{lemma}\label{l:pto} 
Under the assumptions of Proposition~\ref{p:massest}, 
$$\|[\theta(H),\cchi_R(A)]f_n\|= O_R(R^{s-1}),\, \|[\theta(H),C_R]f_n\|=
O_R(R^{s-1}),$$
$$\|[\theta(H),\langle A\rangle^{-s}\tilde\cchi_R(A)] f_n\| = O_R(R^{-1}).$$
\end{lemma}  
\proof By Corollary \ref{co:step0}, the families 
$(\tilde\cchi_R)_{R>1}$ and $(\psi_R^{1/2})_{R>1}$ are bounded in $\Sc^0$, while 
the family $(\sigma _R:t\mapsto \langle t\rangle^{-s} \tilde\cchi_R(t))_{R>1}$ is 
bounded in $\Sc^{-s}$. Furthermore, $]-R, R[$ does not intersect the
supports of $\tilde\cchi_R$, $\psi_R^{1/2}$, and $\sigma_R$. By Lemma~\ref{l:est3}, 
\begin{eqnarray*}  
\|[\theta(H),\cchi_R(A)]\langle A\rangle^{s}\|=
\|[\theta(H),\tilde\cchi_R(A)]\langle A\rangle^{s}\|=
O(R^{s-1}), \\
\|[\theta(H),C_R]\langle A\rangle^{s}\| =O(R^{s-1}), \, \|[\theta(H),\langle A\rangle^{-s}
\tilde\cchi_R(A)] \langle A\rangle^{s}\|= O(R^{-1}). 
\end{eqnarray*}  
Using the boundness of $(\|\langle A\rangle^{-s}f_n\|)_n$ (cf. 
Definition~\ref{suite-speciale}), this yields the results. \qed

\begin{lemma}\label{l:mino-commut}
The inequality \eqref{eq:1} holds true. 
\end{lemma}  
\proof Applying \eqref{eq:mourrer} (with $K=0$) to the 
$\tilde\cchi_R(A)\langle A\rangle^{-s} f_n$, 
\begin{eqnarray*} 
\langle \theta(H) \tilde\cchi_R(A)\langle A\rangle^{-s} f_n, [H_\tau,iA] 
\theta(H) \tilde\cchi_R(A)\langle A\rangle^{-s} f_n\rangle\\ 
&\hspace*{-9cm}\geq&\hspace*{-4cm} c 
\|\theta(H)\tilde\cchi_R(A)\langle A\rangle^{-s} f_n\|^2.    
\end{eqnarray*}  
Recall that $\theta(H)f_n=f_n$. By Lemma~\ref{l:pto}, 
\begin{eqnarray*} 
\left\{\begin{array}{l} 
|\langle [\theta(H),\tilde\cchi_R(A)\langle A\rangle^{-s}] f_n, [H_\tau,iA] 
 \tilde\cchi_R(A)\langle A\rangle^{-s} f_n\rangle |,\\ 
|\langle [\theta(H),\tilde\cchi_R(A)\langle A\rangle^{-s}] f_n, [H_\tau,iA] 
\theta(H)\tilde\cchi_R(A)\langle A\rangle^{-s} f_n\rangle |,\\
|\langle [\theta(H),\tilde\cchi_R(A)\langle A\rangle^{-s}] f_n, 
 \tilde\cchi_R(A)\langle A\rangle^{-s} f_n\rangle |,\\ 
|\langle [\theta(H),\tilde\cchi_R(A)\langle A\rangle^{-s}] f_n, 
\theta(H)\tilde\cchi_R(A)\langle A\rangle^{-s} f_n\rangle |
\end{array}\right. 
\end{eqnarray*}  
are bounded above by $O_R(R^{-1})\cdot \|\tilde\cchi_R(A) \langle A\rangle^{-s} f_n\|$.
Therefore, 
\begin{eqnarray*} 
\langle \tilde\cchi_R(A)\langle A\rangle^{-s} f_n, [H_\tau,iA] 
\tilde\cchi_R(A)\langle A\rangle^{-s} f_n\rangle\\  
&\hspace*{-12cm}\geq&\hspace*{-6cm} c 
\|\tilde\cchi_R(A)\langle A\rangle^{-s} f_n\|^2 +O_R(R^{-1})\cdot 
\|\tilde\cchi_R(A)\langle A\rangle^{-s} f_n\| 
\end{eqnarray*}  
By Lemma \ref{l:gdo2}, for $c'= (2s-1)^{-1}c/2$, \eqref{eq:1} holds true for $R\geq R_1>2$, 
if $R_1$ is large enough. \qed

\begin{lemma}\label{l:mino-C-R}
The inequality \eqref{eq:mino-C-R} holds true.  
\end{lemma}  
\proof From \eqref{eq:mourrer} (with $K=0$) applied to the $C_R f_n$, where $C_R$ is defined in 
\eqref{eq:cr}, we derive that 
\begin{eqnarray*} 
\langle \theta(H) C_R f_n, [H_\tau, iA] \theta(H)  C_R f_n \rangle  \geq 0. 
\end{eqnarray*}  
Thanks to \eqref{eq:cr} and to the Lemmata~\ref{l:pto} and~\ref{l:propC}, 
\begin{eqnarray*} 
\left\{\begin{array}{l} 
|\langle [\theta(H),C_R]  f_n, [H_\tau , iA]C_R\langle A\rangle^{s}\tilde\cchi_{R/2}(A)\langle A\rangle^{-s}
f_n\rangle |\, ,\\ 
|\langle [\theta(H),C_R] f_n, [H_\tau, iA] \theta(H) C_R \langle A\rangle^{s} \tilde\cchi_{R/2}(A) 
\langle A\rangle^{-s}  f_n \rangle | 
\end{array}\right.   
\end{eqnarray*}  
are bounded by $O_R(R^{2s-1}) \cdot\|\tilde\cchi_{R/2}(A) \langle A\rangle^{-s} f_n\|$, 
yielding \eqref{eq:mino-C-R}. \qed

\subsection{Absence of mass.}  

The aim of this part is to show  
\begin{proposition}\label{p:mass2} 
Under the hypotheses of Proposition~\ref{p:massest} with $1/2\l s\l 
2/3$, 
\begin{eqnarray}\label{eq:mass2} 
\lim_{R\rightarrow\infty}\limsup_{n \rightarrow \infty}^{}\|\cchi _R(A)f_n \|= 0.  
\end{eqnarray}   
\end{proposition}  
\proof Applying \eqref{eq:mourrer} (with $K=0$) to the $\cchi _R(A)f_n$, 
\begin{eqnarray*} 
\langle \cchi_R(A) f_n, \theta (H)[H_\tau, iA ]\theta (H)\cchi_R(A) f_n\rangle &\geq & 
c\| \theta (H)\cchi_R(A) f_n\| ^2.   
\end{eqnarray*}  
By Lemma~\ref{l:pto}, 
\begin{eqnarray*} 
\left\{\begin{array}{l} 
|\langle [\theta(H),\cchi_R(A)] f_n, [H_\tau,iA] 
 \cchi_R(A) f_n\rangle |,\\ 
|\langle [\theta(H),\cchi_R(A)] f_n, [H_\tau,iA] 
\theta(H)\cchi_R(A)f_n\rangle |,\\
|\langle [\theta(H),\cchi_R(A)] f_n, 
 \cchi_R(A)f_n\rangle |,\\ 
|\langle [\theta(H),\cchi_R(A)] f_n, 
\theta(H)\cchi_R(A)f_n\rangle |
\end{array}\right. 
\end{eqnarray*}  
are bounded by $O_R(R^{s-1})\cdot \|\cchi_R(A)f_n\|$. Therefore, 
\begin{eqnarray*} 
\langle \cchi_R(A) f_n, [H_\tau, iA ]\cchi_R(A) f_n\rangle \geq 
c\| \cchi_R(A) f_n\| ^2 + O_R(R^{s-1})\cdot \|\cchi_R(A)f_n\| .   
\end{eqnarray*}  
Since $s\l 2/3$, we can find $\beta>0$ (see Lemma~\ref{l:im} below) such that 
\begin{eqnarray} \label{eq:im}
|\langle [H_\tau, \cchi_R(A)]f_n, iA\cchi_R(A) f_n\rangle|\leq  O_R(R^{-\beta})
\|\cchi_R(A) f_n\| . 
\end{eqnarray}  
This yields 
\begin{eqnarray*} 
\langle f_n, [H_\tau, iA \cchi _R^2(A) ]f_n\rangle \geq 
c\| \cchi_R(A) f_n\| ^2 + o_R(1)\cdot \|\cchi_R(A)f_n\| .   
\end{eqnarray*}  
Now, we combine \eqref{eq:carre} and Proposition~\ref{p:virial} to arrive at 
\begin{eqnarray*} 
\limsup _{n\to \infty}\| \cchi_R(A) f_n\| = o_R(1).   \qed
\end{eqnarray*}  

To complete the proof of Proposition~\ref{p:mass2}, we show 

\begin{lemma}\label{l:im} 
Under the assumptions of Proposition~\ref{p:mass2}, there exists $\beta >0$ such that \eqref{eq:im} holds true. 
\end{lemma}  
\proof We decompose $\langle [H_\tau, \cchi_R(A)]f_n, iA\cchi_R(A)f_n\rangle$ as 
\begin{eqnarray}
 \langle [H_\tau, \tilde\cchi_R(A)]
\tilde\cchi _{R/2}(A) f_n, iA\cchi_R(A) f_n\rangle\label{eq:comm-tilde-chi}\\
+ \langle [H_\tau, \tilde\cchi_R(A)]
\cchi _{R/2}(A) f_n, iA\cchi_R(A) f_n\rangle \label{eq:comm-chi}
\end{eqnarray} 
Since $(\tilde\cchi_R)_R$ is bounded in $\Sc^{0}$ (cf. Corollary \ref{co:step0}) and since the support 
of $\tilde\cchi_R$ does not intersect $]-R,R[$, Lemma~\ref{l:est3} for $k=1$ ensures that 
$A[H_\tau, \tilde\cchi_R(A)]\langle A \rangle^{s}$ is bounded and its norm is $O(R^s)$. 
Since $s<2/3$, we can find $\alpha\in ]-s,2s-2[$. This implies, using Corollary \ref{p:mass}, 
that the absolute value of \eqref{eq:comm-tilde-chi} is bounded by $O_R(R^{s+\alpha})\cdot \| \cchi_R(A) f_n\|$, 
with $s+\alpha<0$. By Proposition~\ref{p:regu} with $k=2$, 
\begin{eqnarray*} 
[H_\tau, \tilde\cchi_R(A)]\cchi _{R/2}(A)=[H_\tau,A]\tilde\cchi_R'(A)\cchi _{R/2}(A)+I_2\cchi
_{R/2}(A)=I_2\cchi _{R/2}(A)
\end{eqnarray*}  
since $\supp \tilde\cchi_R'\cap \supp \tilde\cchi _{R/2}=\emptyset$. Lemma~\ref{l:est3} for $k=2$ 
implies that $AI_2\cchi_{R/2}(A)\langle A \rangle^{s}$ is bounded and its norm is $O(R^{s-1})$. 
In particular, the absolute value of \eqref{eq:comm-chi} is bounded by $O(R^{s-1})\cdot \| \cchi_R(A)f_n\|$. \qed

\section{The LAP for the reduced resolvent.} \label{etude-tal-reduit}

\subsection{Motivation}\label{s:motiv}

An interesting consequence of the LAP \eqref{eq:lap} is the following
propagation estimate (cf. Kato's local smoothness in \cite{abg,jmp,rs4}): there exists $C>0$ such 
that, for all $f\in \Hr$, 
\begin{eqnarray}\label{eq:prop}
\int_{-\infty}^{\infty}\|\langle A \rangle^{-s} e^{itH}E_{\Ic}(H) f\|^2dt\leq
C\|f\|^2.
\end{eqnarray} 
For $E_{\Ic}(H)f\neq 0$, the state $e^{itH}E_{\Ic}(H)f$ must move to ``regions where $|A|$ is large'' 
when $t\to -\infty$ and $t\to +\infty$, since the integral converges. If $Hf=\lambda f$ with 
$\lambda\in \Ic$ and $f\neq 0$, then $e^{itH}E_{\Ic}(H)f=e^{it\lambda}f$, $\|\langle A \rangle^{-s}
e^{itH}E_{\Ic}(H) f\|=\|\langle A \rangle^{-s}f\|$, and the integral in 
\eqref{eq:prop} diverges. Therefore, the LAP cannot hold true near an eigenvalue. However it is 
interesting to find out whether the estimate \eqref{eq:prop} holds true for nonzero states $E_{\Ic}(H)f$ 
which are orthogonal to the eigenvectors associated to eigenvalues in $\Ic$, i.e. nonzero states 
$P^\perp E_{\Ic}(H)f$. Now the reduced LAP \eqref{eq:tal-reduit} on $\Ic'$ with $\Ic\subset 
\mathring{\Ic'}$ implies that 
\begin{equation*}
 \sup_{{\rm Re}z\in \Ic , {\rm Im}z\neq 0}\|\langle A \rangle^{-s}(H-z)^{-1} P^\perp E_{\Ic '}(H)
\langle A\rangle^{-s}\|<\infty 
\end{equation*} 
since $(H-z)^{-1}E_{\R\setminus\Ic '}(H)$ is uniformly bounded, yielding \eqref{eq:prop} 
with $f$ replaced by $P^\perp f$ by Kato's local smoothness (cf. \cite{abg,rs4}). 
Theorem~\ref{th:mourre2} gives a situation where the latter estimate holds true. 

\subsection{Eigenvectors' regularity.}\label{s:regu-vect-pr}

Here we extend the result of \cite{c} on the regularity w.r.t.\ $A$ of eigenvectors of $H$. 
\begin{proposition}\label{p:mourre-regu}
Let $\Ic$ be a bounded, open interval that is included in the continuous spectrum of $H$. Let 
$\Ic ''$ be an open interval such $\overline{\Ic}\subset\Ic ''$. Let $H\in \Cc^k_{\Ic ''}(A)$ 
with integer $k\geq 2$. Assume that the Mourre estimate \eqref{eq:mourre2} holds true on 
$\Ic$. Then, for any eigenvector $f$ of $H$ such that $E_{\Ic}(H)f=f$, $f\in \Dc (A^{k-2})$. 
\end{proposition}  
\proof Let $\tau\in\Cc^\infty_c(\Ic '')$ such that $\tau=1$ near $\Ic$. Take $H_\tau$ as in \eqref{eq:def-H-tau}.  Let $f$ be an eigenvector of $H$ such that $E_{\Ic}(H)f=f$. It is also an 
eigenvector of $H_\tau$ with same eigenvalue. As in the proof of Lemma~\ref{l:mourrer}, 
we may replace $H$ par $H_\tau$ in the commutator in \eqref{eq:mourre2}. Now, 
we can follow the proof in \cite{c}, since $H_\tau\in \Cc^k (A)$. \qed

\begin{rem}\label{r:regu-vect-propres}
Proposition~\ref{p:mourre-regu} extends the result in \cite{c} since we 
only assume the ``local" regularity $H\in \Cc^k_{\Ic ''}(A)$. 
\end{rem}

\subsection{Proof of Theorem~\ref{th:mourre2}}\label{s:preuve-TAL-r}
As in the proof of Theorem~\ref{th:mourre} in Subsection~\ref{local}, $H\in \Cc^2_{\Ic ''}(A)$, 
for any open interval $\Ic ''$ with $\overline{\Ic}\subset\Ic ''$. Let $\theta\in\Cc^\infty_c(\Ic)$. 
In particular, $\theta (H)\in\Cc^2(A)$. By Remark~\ref{r:viriel+mourre}, \eqref{eq:mourre2} 
implies that $PE_\Ic(H)$ is a finite dimensional projection. By Proposition~\ref{p:rang1}, 
the hypothesis ${\rm Ran}PE_\Ic(H)\subset \Dc (A^2)$ implies that $PE_\Ic(H)\in \Cc^2(A)$
Since $\theta (H)P=\theta (H)PE_\Ic(H)$, $\theta(H)P\in \Cc^2(A)$. Let 
$\tau\in \Cc_c^\infty(\Ic '')$ such that $\tau =1$ near $\Ic$. Let $\Ic'$ be closed 
with $\Ic'\subset\mathring{\Ic}$ and $\theta\in \Cc_c^\infty(\Ic)$ with $\theta=1$ near 
$\Ic'$. By Lemma~\ref{l:mourrer}, \eqref{eq:mourrer} holds true (with $H_\tau$ defined in
\eqref{eq:def-H-tau}). Thus 
\begin{eqnarray}
P^\perp \theta(H)[H_\tau, iA]  \theta(H)P^\perp &\geq &c
(\theta(H)P^\perp)^2\nonumber\\
&&+\theta(H)P^\perp K\theta(H)P^\perp .\label{eq:m-proj-K} 
\end{eqnarray}  
Let $\theta_1\in \Cc_c^\infty(\Ic ')$. Since $P^\perp:1-P$ projects onto the continuous 
spectral subspace of $H$, $\theta _1(H)P^\perp$ converges strongly to $0$ as the support of $\theta_1$ 
shrinks to a point. Since $K$ compact, $\|K\theta _1(H)P^\perp\|$ goes to $0$ in the same limit. 
Multiplying \eqref{eq:m-proj-K} by $\theta _1(H)$ on both sides and taking the support of $\theta _1$ 
small enough inside $\Ic'$, we obtain 
\begin{eqnarray*}
P^\perp \theta _1(H)[H_\tau, iA]  \theta _1(H)P^\perp &\geq &(c/2)
(\theta _1(H)P^\perp)^2.
\end{eqnarray*}  
Around any point of $\Ic'$, we thus can find some infinite interval $\Ic_1\subset\Ic$ 
such that the projected Mourre estimate \eqref{eq:mourremain-proj} below
holds true on $\Ic_1$. By Theorem~\ref{th:mourre2'}, the reduced LAP
holds true on any closed $\Ic_1'$ with $\Ic_1'\subset\mathring{\Ic_1}$. By compacity of $\Ic'$, we get 
the reduced LAP on it. \qed

So the proof of Theorem~\ref{th:mourre2} reduces to the proof of a
local and stronger version of it, namely 

\begin{theorem}\label{th:mourre2'} 
Let $\Ic$ be a bounded, open interval. Let $\Ic ''$ be an open interval such $\overline{\Ic}\subset\Ic ''$. 
Let $H\in \Cc^2_{\Ic ''}(A)$ and assume that, for all $\theta \in \Cc^\infty_c(\Ic)$, 
$\theta (H)P\in \Cc^2(A)$. Let $\tau\in\Cc^\infty_c(\Ic '')$ such that $\tau=1$ near $\Ic$. 
Assume the {\em projected} Mourre estimate 
\begin{eqnarray}\label{eq:mourremain-proj} 
\hspace*{1cm}P^\perp E_\Ic (H)[H_\tau, iA]  E_\Ic(H)P^\perp \geq c E_\Ic(H)P^\perp ,
\mbox { with }  c\g 0,
\end{eqnarray}  
holds  true. Then, for any $s\g 1/2$ and any compact interval $\Ic'$ with $\Ic'\subset\mathring{\Ic}$, 
the reduced LAP \eqref{eq:tal-reduit}, respectively to $(\Ic', s, A)$, holds true for $H_\tau$ and $H$. 
\end{theorem} 

\subsection{Proofs of Theorem~\ref{th:mourre2'}}\label{s:preuve-TAL-r-loc}

We shall give two proofs of Theorem~\ref{th:mourre2'}. The first one is a direct generalization to 
the present context of {\em our} proof of Theorem~\ref{th:main}. The second proof is close to the corresponding 
proof in \cite{cgh} and shows that Theorem~\ref{th:main} actually applies to $HP^\perp$. In Remark~\ref{r:compare}, 
we compare the two proofs. In Remark~\ref{r:s-reduit}, we comment on Sahbani's result (cf. \cite{s}) in this context. 

\proof[First proof of Theorem~\ref{th:mourre2'}] By Remark~\ref{r:TAL-s}, we may assume that $1/2<s<1$. 
Assume the reduced LAP for $H$ false on some $\Ic'\subset\mathring{\Ic}$. Let $\theta\in\Cc_c^\infty(\Ic)$ 
with $\theta =1$ on $\Ic'$. Notice that, since $\theta (H),\theta (H)P\in \Cc^2(A)$, $\theta (H)P^\perp\in \Cc^2(A)$. 
Then, using the proof of Proposition~\ref{p:non-lap} and~\ref{p:easy}, we can find a special
sequence $(f_n,z_n)_n$ for $H_\tau$ with positive mass such that $\theta (H)f_n=f_n=P^\perp f_n$, for all $n$. 
Since 
\begin{eqnarray*}
\langle A\rangle^{-s}f_n=\langle A\rangle^{-s}(H_\tau-z_n)^{-1}P^\perp\langle A\rangle^{-s}
\langle A\rangle^{s}(H_\tau-z_n) f_n, 
\end{eqnarray*}  
the reduced LAP for $H_\tau$ on $\Ic'$ must be false. So it suffices to prove the reduced LAP for 
$H_\tau$ on $\Ic'$. Using Proposition~\ref{p:non-lap} and~\ref{p:easy} in a similar way, we 
can show that the reduced LAP for $H_\tau$ on $\Ic'$ holds true if and only if, for all special 
sequence $(f_n,z_n)_n$ for $H_\tau$ such that $\theta (H)f_n=f_n=P^\perp f_n$, for all $n$, 
its  mass is $0$. Now, we take such a special sequence $(f_n,z_n)_n$. Multiplying \eqref{eq:mourremain-proj} 
on both sides by $\theta (H)$, 
\begin{eqnarray}\label{eq:m-proj-theta} 
P^\perp \theta (H)[H_\tau, iA] \theta (H)P^\perp \geq c (\theta (H)P^\perp)^2.
\end{eqnarray}  
Since $\theta (H)P^\perp\in \Cc^2(A)$, we can follow our proof of
Theorem~\ref{th:main} in Section~\ref{new-proof}, yielding the reduced
LAP for $H_\tau$ on $\Ic'$. \qed

\proof[Second proof of Theorem~\ref{th:mourre2'}] Assume for a while that 
Theorem~\ref{th:mourre2'} holds true if $0\not\in\overline{\Ic}$. Under the assumptions of
Theorem~\ref{th:mourre2'}, we can find some real $\mu$ such that $0\not\in\mu+\overline{\Ic}$. Notice that $H$ 
and $H+\mu$ have the same eigenvalues and eigenvectors and that the eigenvalues of $H$ in $\Ic$ are the 
eigenvalues of $H+\mu$ in $\mu+\Ic$. For any $\varphi:\R\rightarrow \R$ bounded and borelian, 
$\varphi (H)=\varphi ((H+\mu)-\mu)$, a function of $H+\mu$. Thus the assumptions of Theorem~\ref{th:mourre2'} 
are satisfied if $H$ is replaced by $H+\mu$ and $\Ic$ by $\mu+\Ic$, and $0\not\in\mu+\overline{\Ic}$. 
Thus, it suffices to prove it when  $0\not \in \overline{\Ic}$. \\
For any $\theta\in\Cc^\infty(\R\setminus\{0\})$, $\theta (H)P^\perp=\theta (HP^\perp)$ by 
Lemma~\ref{l:sw} below. Thus $HP^\perp\in \Cc^2_{\Ic ''}(A)$. Furthermore, using Lemma~\ref{l:mourrer}, 
we derive from \eqref{eq:mourremain-proj} the estimate, for $\theta\in\Cc_c^\infty(\Ic)$ 
with $\theta =1$ near $\Ic'$,  
\begin{eqnarray}\label{eq:m-p} 
\theta(HP^\perp)[(HP^\perp)_\tau, iA]\theta(HP^\perp)\geq c \theta^2(HP^\perp). 
\end{eqnarray} 
Now, we can apply Theorem~\ref{th:main} to $HP^\perp$ with $\Ic=\theta^{-1}(1)$, 
yielding the LAP for $HP^\perp$ on $\Ic'$. Let $z\in \C$ with $\im(z)\neq 0$. By Feshbach decomposition 
(see \cite{bfs} for instance), 
$(HP^\perp-z)^{-1}P^\perp=(H-z)^{-1}P^\perp$. Let $\re(z)\in \Ic'$ and
$s\in [0;1[ $. Setting $\tilde\theta =1-\theta$, we write  
\begin{eqnarray*}
\langle A\rangle  ^{-s}(H-z)^{-1}P^\perp\langle A\rangle^{-s}=\langle A\rangle  ^{-s}(H-z)^{-1}P^\perp
\tilde\theta (H)\langle A\rangle  ^{-s}\\
+\langle A\rangle  ^{-s}(HP^\perp-z)^{-1}\langle A\rangle  ^{-s}\cdot
\langle A\rangle  ^{s}\theta (H)P^\perp\langle A\rangle  ^{-s}. 
\end{eqnarray*} 
Since $\theta (H)P^\perp\in \Cc^1(A)$, $\langle A\rangle  ^{s}\theta(H)P^\perp\langle A\rangle ^{-s}$ is bounded by
Proposition~\ref{p:regu}. This yieds the reduced LAP \eqref{eq:tal-reduit} for $H$, since 
$(H-z)^{-1}\tilde\theta (H)$ is uniformly bounded for $\re (z)\in \Ic '$. \qed

The second proof of Theorem~\ref{th:mourre2'} uses the following
consequence of the Feshbach decomposition (see \cite{bfs} for
instance). 
\begin{lemma}\label{l:sw} 
For all $\varphi\in\Cc_c^\infty(\R\setminus\{0\})$, $\varphi(HP^\perp)P=0$ and $\varphi(HP^\perp)=
\varphi(H)P^\perp$.   
\end{lemma}  
\proof  
Let $z\in \C$ with $\im(z)\neq 0$. By Feshbach decomposition, $(HP^\perp-z)^{-1}P^\perp=(H-z)^{-1}P^\perp$. 
Using \eqref{eq:int}, $\varphi(HP^\perp)P^\perp=\varphi(H)P^\perp$. Since ${\rm Ran}P$ is contained in 
the kernel of $HP^\perp$, $\varphi(HP^\perp)P=\varphi(0)P=0$, by assumption on $\varphi$. Finally, 
$\varphi(H)P^\perp=\varphi(HP^\perp)-0$.\qed 

\begin{rem}\label{r:compare}
In the second proof, the idea is to replace $H$ par $HP^\perp$. Since we push that way the 
eigenvectors of $H$ leaving in ${\rm Ran}E_\Ic(H)$ in the kernel of $HP^\perp$, they are 
no longer an obstacle to the strict Mourre estimate on $\Ic$, if $0\not\in\overline{\Ic}$. The main difference between the two previous proofs is probably the use of energy translation for $H$ in 
the second one to avoid the case where $0\in \overline{\Ic}$. 
\end{rem} 
\begin{rem}\label{r:s-reduit}
In the second proof, $HP^\perp \in \Cc ^2_{\Ic "}(A)$ and \eqref{eq:m-p} can be written as 
a Mourre estimate for $HP^\perp $. Under the assumptions of Theorem~\ref{th:mourre2'}, 
Sahbani's result in \cite{s} applies and the boundary values of the reduced resolvent 
have some H\"older continuity. As shown in the proof of Theorem~\ref{th:mourre2}, 
the assumption ``$\theta (H)P\in \Cc ^2(A)$" is satisfied if the Mourre estimate \eqref{eq:mourre2} 
holds true on $\Ic$, included in the continuous spectrum of $H$, and if $H\in \Cc ^4_{\Ic "}(A)$, by Proposition~\ref{p:mourre-regu}. 
\end{rem} 

\subsection{An artificial but instructive example. }\label{artificiel}

In this section, we construct an example of operators $H$ and $A$, for which Theorems~\ref{th:mourre2} 
and~\ref{th:mourre2'} apply but the Mourre estimate \eqref{eq:mourre2} cannot be true. In particular, 
Theorems~\ref{th:mourre} and~\ref{th:main} do not apply to this example. Our contruction is quite 
artificial but our operators $H$ and $A$ presents some structural similarity with operators in \cite{dj}. 

Let $\Hr_0,\Hr_1$ be infinite dimensional complex Hilbert spaces. Let $H_0$ and $A_0$ be 
self-adjoint operators in $\Hr_0$ such that $H_0$ is bounded, $H_0\in \Cc^2(A_0)$, 
and such that the strict Mourre estimate \eqref{eq:mourre2} with $K=0$ holds true for $H_0$
and $A_0$ on some bounded, infinite interval $\Ic$. 
For instance, we can take suitably a bounded, infinite interval $\Ic$ included in $]0;+\infty[$, 
$\Hr_0={\rm L}^2(\R^d)$, $H_0$ a smooth, increasing, and bounded 
function of the Laplacian on $\R^d$, and $A_0$ the generator of dilation in $\R^d$ (cf. \cite{abg,m}). 
Let $A_1$ be self-adjoint operator in $\Hr_1$. 
Let $(g_n)_n$ be a bounded sequence in $\Dc(A_1^2)$ of
independent vectors such that it is bounded for the graph norm of $A_1^2$. 
Let $(\alpha _n)_n\in \ell^1$, a sequence of nonzero reals. 
The serie $(\sum _{n\geq 0}\alpha_n |g_n\rangle \langle g_n|)_n$
converge absolutely in the Banach space of bounded operators on
$\Hr_1$. Let $C$ be its sum. It is a self-adjoint, compact operator of infinite rank. 
By Proposition~\ref{p:rang1}, each $\alpha_n |g_n\rangle \langle g_n|\in \Cc^2(A_1)$ 
and $(\sum _{n\geq 0}\alpha_n [|g_n\rangle \langle g_n|,A_1])_n$ converges absolutely, 
since $(\|g_n\|)_n$ and $(\|A_1g_n\|)_n$ are bounded. By Proposition~\ref{p:fermeture}, 
$C\in \Cc^1(A_1)$ and $[C,A_1]=\sum _{n=0}^\infty\alpha_n 
[|g_n\rangle \langle g_n|,A_1]$. 
Applying this argument again, this implies that $C\in \Cc^2(A_1)$. Let $\lambda\in\mathring{\Ic}$. 
We can choose $(\alpha _n)_n$ such that $[\lambda -\|C\|;\lambda +\|C\|]\subset \Ic$. 
Let $H_1=\lambda+C$. Let 
$H$ be the bounded self-adjoint operator acting in $\Hr:=\Hr_0\oplus \Hr_1$ by 
$H_0\oplus H_1$. Let $A$ be the self-adjoint operator acting in $\Hr$ by
$A_0\oplus A_1$. Since $[H,iA]=[H_0,iA_0]\oplus [C;iA_1]$ as form on
$\Dc(A)\times \Dc(A)$, the regularity of $H_0$ w.r.t.\ $A_0$ and the regularity of 
$C$ w.r.t.\ $A_1$ imply that $H\in \Cc^2(A)$. Since ${\rm Ran}C$ is infinite dimensional 
and the spectrum of $H_1$ is contained in $\Ic$, the point spectrum of $H$ in $\Ic$ is infinite 
therefore the Mourre estimate \eqref{eq:mourre2} cannot hold true on $\Ic$
by Remark~\ref{r:viriel+mourre}. Since the strict Mourre estimate for
$H_0$ holds true on $\Ic$, $H_0$ has no eigenvalue in $\Ic$ by
Remark~\ref{r:viriel+mourre}. 
Let $P$ be the orthogonal projection onto the pure point spectral
subspace of $H$. By the previous properties, 
$P^\perp[H,iA]P^\perp=P^\perp([H_0,iA_0]\oplus 0)P^\perp$. 
Thus the strict Mourre estimate for $H_0$ on $\Ic$ implies the strict, projected 
Mourre estimate \eqref{eq:mourremain-proj} for $H$ on $\Ic$.

\appendix 
\renewcommand{\theequation}{\thesection .\arabic{equation}} 
\section{Symbolic calculus. }\label{calcul-symb}
\setcounter{equation}{0}

In this section, we recall well known facts on symbolic calculus and 
almost analytic extensions (see \cite{dg}[Appendix C]). We also show that some sequences of
functions used in the main text are bounded in some symbol class. 
 
For $\rho\in\R$, let $\Sc^\rho$ be the class of function $\varphi\in\Cc^\infty(\R;\C)$ such that  
\begin{eqnarray}\label{eq:regu} 
\forall k\in\N, \quad C_k(\varphi) :=\sup _{t\in\R}\, \langle t\rangle^{-\rho+k}|\partial_t^k
 \varphi(t)|<\infty .     
\end{eqnarray} 
We also write $\varphi^{(k)}$ for $ \partial_t^k\varphi$. 
Equiped with the semi-norms defined by (\ref{eq:regu}), $\Sc^\rho$ is
a Fr\'echet space. Leibniz' formula implies the continuous embedding: 
\begin{eqnarray}\label{eq:mult} 
 \Sc^\rho\cdot  \Sc^{\rho'} \subset  \Sc^{\rho+\rho'}. 
\end{eqnarray}  
 
For the functional calculus of the operator $A$ (see (\ref{eq:int} )),
we shall use the following result in \cite{dg} on almost analytic extension. 
 
\begin{lemma}\label{l:dg} 
Let $\varphi\in\Sc^\rho$ with $\rho\in\R$. For all   
$l\in \N$, there is a smooth function  $\varphi^\C:\C \rightarrow \C$, call an 
\emph{almost analytic extension} of $\varphi$, such that: 
\begin{eqnarray} 
\label{eq:dg1} \varphi^\C|_{\R}=\varphi,\quad &&\big|\frac{\partial 
  \varphi^\C}{\partial  \overline{z}}(z) \big|\leq c_1 \langle \re(z) 
  \rangle^{\rho-1 -l} |\im(z)|^l\\\label{eq:dg2} 
&& \supp \varphi^\C \subset\{x+iy\mid |y|\leq c_2 \langle 
  x\rangle\},\\\label{eq:dg3} && \varphi^\C(x+iy)= 0, \mbox{ if } 
  x\not\in\supp\varphi .  
\end{eqnarray}  
for constants $c_1$, $c_2$ depending on the semi-norms 
\eqref{eq:regu} of $\varphi$ in $\Sc^\rho$.   
\end{lemma}  
  
The function $\cchi_R$, given by \eqref{eq:chi}, belongs to $\Sc^\rho$, 
for any $\rho$ and any $R$. But we need to know that the family 
$(\chi_R)_{R\geq 1}$ is bounded in some $\Sc^\rho$. 

\begin{lemma}\label{l:uni} 
Let $\tau\in \Cc^\infty(\R;\R)$ such that $\tau'\in \Cc_c^\infty(\R^\ast;\R)$. Then the family
$(\tau_R)_{R\g 1}$, with $\tau_R(x):=\tau (x/R)$, is bounded in $\Sc^{0}$. 
\end{lemma}  
\proof Let $k\in\N$. The semi-norm $C_k(\tau)$ (cf. \eqref{eq:regu}) is bounded above
by the $(\sup {\rm supp}\tau ^{(k)})^k$ times the $L^\infty$-norm of
$\tau ^{(k)}$. For $R\g 1$ and $t\in\R$, 
\begin{eqnarray}\label{eq:deriv-tau-R} 
|t|^k\cdot |(\tau_R)^{(k)}(t)|=(|t|/R)^k\cdot |\tau ^{(k)}(t/R)|\leq
C_k(\tau).
\end{eqnarray}  
Thus $(\tau_R)_{R\g 1}$ is bounded in $\Sc^{0}$. \qed 

Concerning the functions defined in \eqref{eq:chi}, \eqref{eq:chi2}
and just after \eqref{eq:chi2}, we have the 

\begin{corollary}\label{co:step0}  
Lemma~\ref{l:uni} applies to $\tau=\cchi,\tilde\cchi,\tilde\cchi
_\sigma, (\tilde\cchi _\sigma)^{1/2}$, $(\sigma \tilde\cchi
_\sigma')^{1/2}$, for $\sigma\in \{-;+\}$, and also to their 
derivatives. 
\end{corollary}  

We now focus on the functions $\psi_{\sigma,R}$, defined in \eqref{eq:positive}. 
 
\begin{lemma}\label{l:step} 
The family $(\psi_{\sigma,R}^{1/2})_{R\g 1}$ is bounded in $\Sc^{0}$.    
\end{lemma}  
\proof 
By \eqref{eq:positive},   
\begin{eqnarray*} 
\psi_{\sigma,R}^{1/2}(x)=(\sigma\cchi'_{\sigma})^{1/2}(x/R)
(h -\sigma x \langle x\rangle^{-2s})^{1/2}\cchi_{\sigma}^{1/2}(x/R), 
\end{eqnarray*}  
for all $x\in \R$ and all $R\g 1$. By definition of $h$ (cf. \eqref{eq:h}), 
$x\mapsto (h -\sigma x \langle x\rangle^{-2s})^{1/2}$ belongs to $\Sc^{0}$. 
Now the result follows from Corollary \ref{co:step0} and \eqref{eq:mult}.\qed

\section{Commutator expansions. } \label{dev-commut}
\setcounter{equation}{0}

In this section, we recall Helffer-Sj\"{o}strand's functional calculus
(cf. \cite{hs,dg}) and commutator expansions (cf. \cite{dg}). 

Let $\rho\l 0$ and $\varphi\in \Sc^{\rho}$. The bounded operator $\varphi(A)$ 
can be recover by Helffer-Sj\"{o}strand's formula: 
\begin{eqnarray}\label{eq:int} 
\varphi(A) = \frac{i}{2\pi}\int_\C\frac{\partial \varphi^\C }{\partial 
\overline{z}}(z-A)^{-1}dz\wedge d\overline{z}, 
\end{eqnarray}  
where the integral exists in the norm topology, by \eqref{eq:dg1} 
with $l=1$. This can be extended as shown in 
\begin{lemma}\label{l:ext-hs} 
Let $k\in \N$, $\rho <k$ and $\varphi\in \Sc^{\rho}$. Strongly in $\Dc(\langle A\rangle^{k})$,
\eqref{eq:int} holds true.  
\end{lemma}  
\proof Let $f\in \Dc(\langle A\rangle^{k})\subset \Dc(\varphi(A))$ 
and $\cchi_R$ be like in \eqref{eq:chi}, then 
\begin{eqnarray}\label{eq:chiRf}
\hspace*{1cm}\varphi(A)\cchi_R(A)f = \frac{i}{2\pi}\int_\C\frac{\partial 
(\varphi _k\cchi_R)^\C }{\partial \overline{z}}(z-A)^{-1}\langle A\rangle^{k}f\, dz\wedge d\overline{z}, 
\end{eqnarray} 
where $\varphi _k(t):=\varphi (t)\langle t\rangle ^{-k}$ belongs to $\Sc^{\rho -k}$. 
By Lemma~\ref{l:uni} and \eqref{eq:mult}, $(\varphi _k\cchi_R)_R$ is
bounded in $\Sc^{\rho -k}$. Since $\rho -k<0$, the result follows from \eqref{eq:dg1} with $l=1$
and the dominated convergence theorem. \qed 

Notice that, for some $c>0$ and $s\in [0;1[$, there exists
some $C>0$ such that, for all $z=x+iy\in\{a+ib\mid 0<|b|\leq c\langle a\rangle \}$ 
(like in \eqref{eq:dg2}), 
\begin{eqnarray}\label{eq:majoA} 
\big\| \langle A\rangle^s (A-z)^{-1}\big\|\leq C \langle x \rangle^{s}\cdot |y|^{-1}. 
\end{eqnarray}   

Next we come to a commutator expansion. 
 
\begin{proposition}\label{p:regu} 
Let $k\in \N^\ast$ and $B$ be a self-adjoint and bounded operator in $\Cc^k(A)$. Let
$\rho\l k$ and $\varphi\in \Sc^{\rho}$. In the sense of forms on 
$\Dc(\langle A\rangle^{k-1})\times \Dc(\langle A\rangle^{k-1})$: 
\begin{eqnarray}\label{eq:egalite}
[\varphi(A), B] = \sum_{j=1}^{k-1} \frac{1}{j!} 
\varphi^{(j)}(A)\ad_A^j(B)\\ \label{eq:reste22}    
+ \frac{i}{2\pi}\int_\C\frac{\partial\varphi^\C }{\partial 
\overline{z}}(z-A)^{-k} \ad_A^k(B) (z-A)^{-1} dz\wedge d\overline{z}. 
\end{eqnarray} 
In particular, if $\rho\l 1$, then $B\in\Cc^1(\varphi(A))$. 
\end{proposition} 
\proof  

Thanks to Lemma~\ref{l:ext-hs}, we can write, as form on 
$\Dc(\langle A\rangle^{k})\times \Dc(\langle A\rangle^{k})$: 
\begin{eqnarray*}  
\hspace*{1cm} 
[\varphi(A), B] &=& \frac{i}{2\pi}\int_\C\frac{\partial\varphi^\C }{\partial 
\overline{z}}(z-A)^{-1} \ad_A(B) (z-A)^{-1} dz\wedge 
d\overline{z}\\ 
&=&\sum_{j=1}^{k-1} \frac{i}{2\pi}\int_\C\frac{\partial\varphi^\C 
}{\partial\overline{z}}  (z-A)^{-j-1}\ad_A^{j}(B) dz\wedge 
d\overline{z}\\ 
&+& \frac{i}{2\pi}\int_\C\frac{\partial\varphi^\C }{\partial 
\overline{z}}(z-A)^{-k} \ad_A^k(B) (z-A)^{-1} dz\wedge d\overline{z}. 
\end{eqnarray*} 
This yields \eqref{eq:egalite} on $\Dc(\langle A\rangle^{k})\times \Dc(\langle A\rangle^{k})$. 
Since $B\in \Cc^k(A)$, the commutators $\ad_A^{j}(B)$, for $1\leq j\leq k$, are bounded. 
Now, as in the proof of Lemma~\ref{l:ext-hs}, we see that this form extends to a bounded form on 
$\Dc(\langle A\rangle^{k-1})\times \Dc(\langle A\rangle^{k-1})$ since 
the $\varphi^{(j)}$ belong to $\Sc^{\rho -1}$.  \qed

The rest of the previous expansion is estimated in 
\begin{lemma}\label{l:est3} 
Let $B\in\Cc^k(A)$ self-adjoint and bounded. Let $\varphi\in\Sc^\rho$, with $\rho\l k$. Let $I_k(\varphi)$ the 
rest of the development of order $k$ \eqref{eq:egalite} of $[\varphi(A), B]$, namely \eqref{eq:reste22}. Let 
$s, s'\geq 0$ such that $s'<1$, $s<k$, and $\rho+s+s'<k$. Then $\langle A \rangle^{s} 
I_k(\varphi)\langle A \rangle^{s'}$ is bounded and it is uniformly bounded when
$\varphi$ stays in a bounded subset of $\Sc^\rho$. In particular, $I_k(\varphi)$ is a bounded 
operator. Let $R>0$. If $\varphi$ stays in a bounded subset of $\{\psi
\in \Sc^\rho\mid [-R;R]\cap \supp(\varphi)=\emptyset\}$ then $\langle R\rangle^{k-\rho-s-s'}
\|\langle A \rangle^{s} I_k(\varphi)\langle A \rangle^{s'}\|$ is uniformly bounded. 
\end{lemma}  
\proof 
We will follow ideas from \cite{dg}[Lemma C.3.1]. In this proof, all the 
constants are denoted by $C$, independently of their value. Given a complex 
number $z$, $x$ and $y$ will denote its real and imaginary part,
respectively. Since $B\in\Cc^k(A)$, $\ad^k_A(B)$ is bounded. 
We start with the second assertion. Let $\varphi\in \Sc^\rho$, $R>0$
such that $[-R;R]\cap \supp(\varphi)=\emptyset$. Notice that, 
by \eqref{eq:dg3}, $\varphi^\C(x+iy)= 0$ for $|x|\leq R$. By definition of $I_k$, we
consider \eqref{eq:reste22} and switch to the variable $(x,y)$ by noticing that 
$dz\wedge d\overline{z}=-2i dx\wedge dy$. By \eqref{eq:majoA}, 
\begin{eqnarray*} 
\|\langle A \rangle^{s}I_k(\varphi)\langle A \rangle^{s'}\|\leq \frac{1}{\pi} \int 
\big|\frac{\partial\varphi^\C }{\partial \overline{z}}\big|\cdot 
\frac{\langle   x\rangle^{s}}{|y|^k}  \cdot  
\|\ad^k_A(B)\|\cdot \frac{\langle x\rangle^{s'}}{|y|} dx\wedge dy\\ 
\leq  C(\varphi)\int_{|x|\geq R}\int_{|y|\leq c_2\langle x\rangle } \langle x 
\rangle^{\rho+s+s'-1-l}|y|^l |y|^{-k-1} dx\wedge dy,  
\end{eqnarray*}  
for any $l$, by \eqref{eq:dg1}. We choose $l=k+1$. We have,  
\begin{eqnarray*} 
\|\langle A\rangle^{s}I_k(\varphi)\langle A\rangle^{s'} 
\|&\leq& C(\varphi)\int_{|x|\geq R} \langle x \rangle^{\rho+s+s'-k-1}dx 
\\ 
&\leq& C(\varphi)\langle R\rangle^{\rho+s+s'-k}.   
\end{eqnarray*}  
Since $C(\varphi)$ is bounded when $\varphi$ stays in a bounded subset
of $\Sc^\rho$, this yields the second assertion. For the first one, we can follow the same lines, 
replacing $R$ by $0$ in the integrals, and arrive at the result. \qed

\section{Technical estimates. } \label{esti-techn}
\setcounter{equation}{0}

\begin{lemma}\label{l:gdo2} 
Let $\varepsilon\in ]0,1-s[$ and suppose $B\in \Cc^2(A)$ bounded and 
self-adjoint. Then, for all $f\in\Hr$,  
\begin{eqnarray*} 
 \langle \tilde\cchi_R (A)  f, [B, - A\langle A\rangle^{-2s}  ]\tilde\cchi_R 
 (A)f\rangle &\geq& O_R\big(R^{-\varepsilon})\|\tilde\cchi_R(A) \langle 
 A\rangle^{-s} f\|^2 \\  
&&\hspace*{-5cm}+ (2s-1) \langle \tilde\cchi_R (A)\langle 
 A\rangle^{-s} f, [B, A]  
\tilde\cchi_R(A)\langle A\rangle^{-s} f \rangle .      
\end{eqnarray*}  

\end{lemma}  
\proof  
Let $D= [B, -\langle A\rangle^{-2s}  A]- (2sA^2 \langle 
A\rangle^{-2s-2}-\langle A\rangle^{-2s})[B,iA]$. By Lemma 
\ref{l:est3} for $k=2$, as $t\mapsto \langle t\rangle^{-2s}t\in \Sc^{1-2s}$, 
one has $\langle A\rangle^{s+\varepsilon} D\langle A\rangle^{s}$ bounded for  
$\varepsilon <1$. Then, using the fact that 
$\tilde\cchi_{R/2}(t)=1$  for $t$ in the support of  $\tilde\cchi_R$, 
\begin{eqnarray*} 
\langle \tilde\cchi_R (A)f, D\tilde\cchi_R (A) f\rangle &=&
\langle \langle A\rangle^{-\varepsilon}\tilde\cchi_{R/2} (A)\tilde\cchi_R (A) 
\langle A\rangle^{-s}  
f, \langle A\rangle^{s+\varepsilon}D   \tilde\cchi_R (A) f\rangle\\ 
&\leq &O(R^{-\varepsilon})\cdot \|\tilde\cchi_R (A)\langle 
A\rangle^{-s} f\|^2.    
\end{eqnarray*}  
Since $[B,A]\in \Cc^1(A)$ and $t\mapsto \langle t\rangle^{-s}\in\Sc^{-s}$, 
Lemma \ref{l:est3} gives that $\langle A\rangle^{s+\varepsilon}[\langle 
  A\rangle^{-s}, [B,A]]\langle A\rangle^{s}$ bounded for $\varepsilon 
<1-s$. Using, like above, the  
contribution of $\tilde\cchi_{R/2}(A)$, 
\begin{eqnarray*} 
|\langle \tilde\cchi_R (A)f, ( 2sA^2 \langle 
A\rangle^{-2s-2}-\langle A\rangle^{-2s})[B,iA]\tilde\cchi_R (A)f\rangle 
&&  \\ 
&\hspace*{-16cm}-&\hspace*{-8cm}  
\langle \tilde\cchi_R (A)f, (2sA^2\langle A\rangle^{-2} -1)\langle 
A\rangle^{-s}[B,iA]\langle A\rangle^{-s}\tilde\cchi_R (A)f\rangle |\\  
&\hspace*{-16cm}\leq &\hspace*{-8cm}O(R^{-\varepsilon})\|\langle 
A\rangle^{-s} \cchi_R (A)f\|^{2}.  
\end{eqnarray*}  
To conclude, observe that $\|\tilde\cchi_{R/2}(A)(\id- A^2 \langle 
A\rangle^{-2}) \|= O(R^{-2})$. \qed

\begin{lemma}\label{l:exp} 
Let $B\in\Cc^2(A)$ bounded and self-adjoint. For all $f\in\Hr$,    
\begin{eqnarray*} 
 |\langle f, (\sigma h - A 
\langle A\rangle^{-2s}) \tilde\cchi_{\sigma, R}(A)[B, \tilde\cchi_{\sigma, 
    R}(A)] f \rangle&-&\\  
&\hspace*{-10cm} &\hspace*{-5cm} \langle f, \tilde\cchi_{\sigma, 
  R}'(A)(\sigma h- A \langle A\rangle^{-2s}) \tilde\cchi_{\sigma, R}(A)[B, A] f 
\rangle |\\ 
&\hspace*{-10cm}\leq &\hspace*{-5cm}O(R^{2s-2})\|\tilde\cchi_{R}(A)\langle 
A\rangle^{-s}f\|\cdot \|\langle A\rangle^{-s}f\| . 
\end{eqnarray*} 
\end{lemma}  
\proof By Lemma~\ref{l:est3}, we develop the commutator and denote the rest by $I_2$. 
Its contribution is  
\begin{eqnarray*} 
 \langle \tilde\cchi_{R/2}(A) \langle A\rangle^{-s} f, (\sigma h- A \langle 
 A\rangle^{-2s}) \tilde\cchi_{\sigma, R}(A) \langle A\rangle^s  I_2  
\langle A\rangle^s \langle A\rangle^{-s} f \rangle . 
\end{eqnarray*}  
Note that $\tilde\cchi_{R/2}(A)f$ appears freely thanks to the presence of 
  $\tilde\cchi_{\sigma,  R}(A)$. By Corollary \ref{co:step0}, $(\tilde\cchi_{\sigma, R})_R$ is 
bounded in $\Sc^{0}$. Note also that $[-R,R]$ is not contained in the support of 
$\tilde\cchi_{\sigma, R}$. Then, from  Lemma~\ref{l:est3}, used with $k=2$, we obtain 
  that $\langle A\rangle^s   I_2 \langle A\rangle^s = O(R^{2s-2})$.\qed

\begin{lemma}\label{l:propC} 
For $B\in \Cc^1(A)$ bounded and self-adjoint, 
\begin{enumerate} 
\item $\|C_R \langle A\rangle^{\alpha} \|=O(R^{\alpha})$, for $\alpha\in\R$,   
\item $\|[B,  C_R]\langle A\rangle^{\alpha}\|=O(R^{\alpha-1})$, for $0\leq \alpha < 1$.   
\end{enumerate}   
\end{lemma}  
\proof 
Since $\psi_R(t)=0$ for $|t|\not \in [R,2R]$, the point (1) follows. Since 
$(\psi_R)_R$ is bounded in $\Sc^{0}$ (cf. Lemma~\ref{l:step}) and since $[-R,R]$ is not 
contained in the support of $\psi_R$, we get the point (2) by Lemma~\ref{l:est3}. \qed      
 
\begin{lemma}\label{l:C} 
Let $B\in\Cc^2(A)$ bounded, self-adjoint. For all $f\in\Hr$, 
\begin{eqnarray*} 
  |\langle C_R f, [B,i A] C_R f\rangle &-&   
\langle f, C_R^2[B, iA] f \rangle |\\ 
&&\leq O_R\big(R^{2s-1}) \cdot \|\tilde\cchi_{R/2}(A)\langle A\rangle^{-s}f \| \cdot \|\langle A\rangle^{-s}f \| . 
\end{eqnarray*}  
\end{lemma}  
\proof Given $f\in \Hr$ and using \eqref{eq:cr}, 
\begin{eqnarray*} 
\langle C_R f, [B, iA]C_R f\rangle &=&   
\langle f, C_R^2 [B, iA]  f \rangle  
\\&\hspace*{-8cm}-&\hspace*{-4cm} 
\langle  \tilde\cchi_{R/2}(A) \langle A \rangle^{-s} f, C_R \langle A\rangle^{s}
[C_R,[B, iA]] \langle A\rangle^{s} \cdot \langle A\rangle^{-s} 
f \rangle .       
\end{eqnarray*}  
The last term is estimated above by 
\begin{eqnarray*} 
\|C_R\langle   A\rangle^{s}\|\cdot\|    [C_R,[B, A]] \langle A\rangle^{s}\|
\cdot \|\tilde\cchi_{R/2}(A) \langle  A\rangle^{-s}f\| \cdot \|\langle A\rangle^{-s}f \| .
\end{eqnarray*}  
Now Lemma~\ref{l:propC} gives the result. \qed

\end{document}